\documentclass[MSNbibl,number,citesort,seceqn,dvips]{arxbj}

\aid{0}
\volume{18}
\issue{4}
\pubyear{2012}
\firstpage{1284}
\lastpage{1309}
\doi{10.3150/11-BEJ368}

\makeatletter

\newtheorem{theorem}{Theorem}
\newproclaim{condition}{Condition}
\newtheorem{corollary}{Corollary}
\newproclaim{definition}{Definition}
\newtheorem{lemma}{Lemma}
\newremark{remark}{Remark}

\newcommand{\fracb}[2]{{(#1)}/{#2}}
\makeatother

\begin{document}
\begin{frontmatter}

\title{Continuity and differentiability of regression M functionals}
\runtitle{Continuity and differentiability of regression M functionals}

\begin{aug}
\author[1]{\fnms{Mar\'{i}a V.} \snm{Fasano}\thanksref{1,e1}\ead[label=e1,mark=e1]{virfeather@yahoo.com.ar}},
\author[1]{\fnms{Ricardo A.} \snm{Maronna}\thanksref{1,e2}\ead[label=e2,mark=e2]{rmaronna@retina.ar}},
\author[2]{\fnms{Mariela}~\snm{Sued}\thanksref{2}\ead[label=e3]{msued@dm.uba.ar}}
\and\break
\author[3]{\fnms{V\'{i}ctor~J.}  \snm{Yohai}\corref{}\thanksref{3}\ead[label=e4]{vyohai@dm.uba.ar}}
\runauthor{Fasano, Maronna, Sued and Yohai}

\address[1]{Departamento de Matem\'atica, Facultad de Ciencias Exactas, Universidad Nacional de  La Plata,
Calles 50 y 115, 1900  La Plata,  Argentina. \printead{e1};\printead*{e2}}

\address[2]{Instituto de C\'alculo, Facultad de Ciencias Exactas y Naturales,
Universidad de Buenos Aires, Ciudad Universitaria, Pabell\'{o}n 1, 1426
Buenos Aires, Argentina and CONICET. \printead{e3}}

\address[3]{Departamento de Matem\'atica, Facultad de Ciencias Exactas y Naturales,
Universidad de Buenos Aires, Ciudad Universitaria, Pabell\'{o}n 2, 1426
Buenos Aires, Argentina and CONICET.\\\printead{e4}}
\end{aug}

\received{\smonth{6} \syear{2010}}
\revised{\smonth{2} \syear{2011}}

\begin{abstract}
This paper deals with the Fisher-consistency, weak continuity and
differentiability of estimating functionals corresponding to a class of both
linear and nonlinear regression high breakdown M estimates, which includes  S
and MM estimates. A restricted type of differentiability, called weak
differentiability, is defined, which suffices to prove the asymptotic
normality of estimates based on the functionals. This approach allows to prove
the consistency, asymptotic normality and qualitative robustness of M
estimates under more general conditions than those required in standard
approaches.  In particular, we prove that regression MM-estimates are
asymptotically normal when the observations are $\phi$-mixing.
\end{abstract}


\end{frontmatter}

\section{Introduction}

We consider estimation in the regression model with random predictors
\begin{equation}
y_{i}=g(x_{i},\beta_{0})+u_{i}, \label{modelogral}%
\end{equation}
with data $ (  x_{i},y _{i} )  \in R^{p}\times R$, $i=1,\ldots,n;$
where $\beta_{0}\in B\subseteq R^{q}$ is a vector of unknown parameters,
$g(x,\beta)$ is a known function continuous in $\beta$, and for each $i, x_{i}$
and $u_{i}$ are independent. It is assumed that $\{ (  x_{i}%
,y_{i} )  ,i\geq1\}$ are identically distributed but not necessarily
independent. The well-known fact that the least squares (LS) estimate of
$\beta_{0}$ is sensitive to atypical observations has motivated the
development of robust estimates.

An important class of robust estimators are the \emph{M estimates}. Inside
this class we can distinguish the S estimates introduced by Rousseeuw and
Yohai~\cite{RY} and the MM estimates proposed by Yohai~\cite{YohaiII}. For
linear regression, S estimates may attain the highest possible breakdown
point, and MM estimates may combine the highest possible breakdown point with
a high normal efficiency; see, for example,~\cite{Maronna},
Chapter 5. In the case of nonlinear regression,\ MM estimates may also combine
high breakdown point with a high normal efficiency. In fact, the normal
efficiency of these estimates can be made as close to one as desired, and
Monte Carlo simulations in Fasano~\cite{Fasano} show them to have a highly
robust behavior for some nonlinear models.

In the nonlinear case, Fraiman~\cite{Fraiman} studied bounded influence
estimates for nonlinear regression. Sakata and White~\cite{SakaW} dealt with S
estimates for nonlinear regression models with dependent observations; Vainer
and Kukush~\cite{VaKu} and Liese and Vajda~\cite{LiVa1,LiVa2} dealt
with M estimates with a fixed scale, which therefore lack scale equivariance.
The latter study the $\sqrt{n}$-consistency of M estimates in more general
models, which include linear and nonlinear regression with independent
observations. Stromberg~\cite{StroLMS} proved the weak consistency of the
least median of squares (LMS) estimate$,$ and C\'{i}zek~\cite{CizLTS} dealt
with the consistency and the asymptotic normality of the least trimmed squares
(LTS) estimate under dependency.

Three important qualitative features of an estimate are consistency,
asymptotic normality and qualitative robustness. These properties have been
studied in the literature through specific approaches. Yohai~\cite{YohaiII}
proved these properties for MM estimates in the i.i.d. linear case, and Fasano
\cite{Fasano} proved them in the nonlinear case, both assuming symmetrically
distributed $u_{i}$'s.

In this work, we propose an alternative approach, based on the representation
of the estimates as \emph{functionals} on distributions (Hampel~\cite{Hampel}%
). For a large class of estimates, which includes M estimates, one can define
a functional $T (  G )  $ on the space of data distributions, such
that if $G_{n}$ is the empirical distribution, then $T (  G_{n} )  $
is the estimate, and if $G_{0}$ is the underlying distribution, then $T (
G_{0} )  $ is the parameter that we want to estimate$.$ The weak
continuity of the functional $T$ simplifies the proof of consistency of
$T(G_{n})$ and some suitable forms of differentiability of $T$, as Fr\'{e}chet
or Hadamard differentiability, allow simple proofs of the asymptotic normality
of the estimate under very general conditions. These results hold without the
requirement that $G_{n}$ be the empirical distribution of a sequence of i.i.d.
random variables: if we want to estimate $T(G_{0})$, it suffices that $G_{n}$
converges weakly to $G_{0}$ a.s. The weak continuity of M functionals at a
general statistical model were studied by Clarke~\cite{ClarkeI} and
\cite{ClarkeII}. Fr\'{e}chet differentiability was studied by Boos and
Serfling~\cite{Boos} and Clarke~\cite{ClarkeI}, and Hadamard
differentiability by Fernholz~\cite{Fernholz}. In all of these works, it is
required that the score function used for the M estimate be bounded, and
therefore their results can not be applied to regression M estimates. In this
paper, we prove under very general conditions that the functionals associated
to M estimates of regression are weakly continuous. Besides, since the usual
forms of differentiability, like Fr\'{e}chet or Hadamard differentiability,
require in the case of M estimates the boundedness of the score functions, we
introduce a new concept of differentiability, that we call \emph{weak
differentiability}, which is satisfied by high breakdown M estimates of
regression, for example, by S and MM estimates, and which is adequate to prove the
asymptotic normality of these estimates.

This work is organized as follows: In Section~\ref{SecMain}, we define the
estimates to be considered and in Sections~\ref{SecFish-M},~\ref{SecConti-M-S}
and~\ref{SecDife} we shall respectively deal with the Fisher-consistency,
continuity and differentiability of the functionals corresponding to the
estimates defined above. These results will be shown to imply the consistency,
qualitative robustness and asymptotic normality of the estimates under
assumptions more general than the i.i.d. model and without the requirement of
symmetric errors. In Section~\ref{SectionMM}, we apply the results obtained in
the former sections to MM estimates. Finally, Section~\ref{SecProofs} contains
all proofs.\looseness=-1\vspace*{-1pt}

\section{Definitions of estimates} \label{SecMain}\vspace*{-1pt}

We first define our notation. Henceforth, $\mathrm{E}_{G}[h(z)]$ and
$\mathrm{P}_{G} (  A )  $ will respectively denote the expectation of
$h(z)$ and the probability that $z\in A,$ when $z$ is distributed according to
$G$. If $z$ has distribution $G$, we write $z\sim G$ or $\mathcal{D} (
z )  =G.$ Weak convergence of distributions, convergence in probability
and convergence in distribution of random variables or vectors are denoted by
$G_{n}\rightarrow_{w}G,$ $z_{n}\rightarrow_{p}z$ and $z_{n}\rightarrow_{d}z,$
respectively.  By an abuse of notation, we will write $z_{n}\rightarrow_{d}G$
to denote $\mathcal{D} (  z_{n} )  \rightarrow_{w}G.$ The complement
and the indicator of the set $A$ are denoted by $A^{c}$ and $\mathbf{1}_{A},$
respectively. The scalar product of vectors $a$ and $b$ is denoted by
$a^{\prime}b$, and $R_{+}$ denotes the set of positive real numbers.

To identify $\beta_{0}$, without assuming that the distribution of $u$ is
symmetric around $0$ or that it satisfies a centering condition (such as
e.g. E$_{F_{0}}u=0)$, we assume the following\vspace*{-1pt}

\begin{condition}
\label{CondId} For all $\beta\neq\beta_{0}$ and for all $\alpha$, we have%
\begin{equation}
\mathrm{P} \bigl(  g(x,\beta_{0})=g(x,\beta)+\alpha \bigr)  <1.\vspace*{-1pt} \label{IDCOND}%
\end{equation}
\end{condition}

 Note that when  this condition is not satisfied,  there exist $\beta
\neq\beta_{0}$ and $\alpha$ such that (\ref{modelogral}) also holds with
$\beta$ instead of $\beta_{0}$ and $u_{i}+\alpha$ instead of $u_{i}%
$. Condition~\ref{CondId} requires that in case there is an intercept, it
will be included in the error term $u$ instead  of as a parameter of the
regression function $g(x,\beta)$. For linear regression, we have
$g(x,\beta)=\beta^{\prime}x$ and then this condition means that the vector $x$
is not concentrated on any hyperplane.

Although model (\ref{modelogral}) does not contain an intercept, in order to
 obtain consistent estimates of $\beta_{0},$ our M estimates, besides an
estimate $\widehat{\beta}$ of $\beta_{0}$, will include an additional
additive term $\widehat{\alpha}$. If~the model does contain an intercept, then
$\widehat{\alpha}$ will be a consistent estimate of this parameter under the
centering condition $\mathrm{E}\rho^{\prime}(u/\sigma)=0,$ where $\rho$ is the
loss function of the M estimate and $\sigma$ is the  asymptotic value  of
the  estimate of the error scale that is used to define the M estimate.  If
the model does not contain an intercept, then $\widehat{\alpha}$ can be
ignored.  Let henceforth $\xi=(\beta^{\prime},\alpha)^{\prime}$ with
$\alpha\in R$, and define the function%
\[
\underline{g}(x,{\xi})=g(x,\beta)+\alpha.
\]

\emph{M estimates} are then defined as%
\begin{equation}
\widehat{\xi}_{\mathrm{M}}=\arg\min_{\xi \in B\times R}\sum_{i=1}%
^{n}\rho \biggl(  \frac{y_{i}-\underline{g} (  x_{i},\xi )  }%
{\widehat{\sigma}} \biggr)  , \label{defMestima}%
\end{equation}
where $\widehat{\sigma}$ is a robust residual scale and $\rho$ is a loss
function.\vadjust{\goodbreak}

To define \emph{S estimates}, we need an M scale $S(r)$. Given $r=(r_{1}%
,\ldots,r_{n})^{\prime}$, $S(r)$ is defined as  the solution $\sigma$ of%
\begin{equation}
\frac{1}{n}\sum_{i=1}^{n}\rho_{0} \biggl(  \frac{r_{i}}{\sigma} \biggr)  =\delta,
\label{defMscale}%
\end{equation}
where $\rho_{0}$ is another loss function and the constant $\delta$ regulates
the estimate's robustness.

Then, S estimates of regression are defined by
\begin{equation}
\widehat{\xi}_{\mathrm{S}}=\arg\min_{\xi\in B\times R}S (  r (
\xi )   )  , \label{defSestima}%
\end{equation}
where $r (  \xi )  $ is the residual vector with elements $r_{i}%
(\xi)=y_{i}-$\underline{$g$}$ (  x_{i},\xi )  $.

In particular, we will consider with some detail the subclass of MM estimates.
These estimates  are defined by (\ref{defMestima}) with $\widehat{\sigma}$
obtained from an S estimate, namely
\begin{equation}
\widehat{\sigma}=\min_{\xi\in B\times R}S (  r (  \xi )   )
\label{Sscale}%
\end{equation}
with $ \rho\leq\rho_{0}$. Yohai~\cite{YohaiII} showed  that in the case of
linear regression the  asymptotic breakdown point of MM estimates with
$\delta=0.5$ is $0.5$ if $\mathrm{P}(\beta^{\prime}x_{i}+a=0)=0$ for all
$\beta\neq0$, and that, simultaneously,  it is possible to choose $\rho$ so
that the corresponding MM estimate yields an arbitrarily high efficiency when
the errors are Gaussian.

Now in order to state our results, we must first express the already defined M
and S estimates as functionals. Throughout this article, loss functions will
be ``bounded $\rho$-functions,'' in the
following sense.

\begin{definition}
\label{DefRho}
A \emph{bounded} $\rho$-\emph{function} is a function
$\rho (  t )  $ that is a continuous nondecreasing function of $|t|,$
such that $\rho(0)=0,$ $\rho (  \infty )  =1,$ and $\rho (
v )  <1$ implies that $\rho (  u )  <\rho(v)$ for $|u|<|v|.$
\end{definition}

Then, in the rest of the paper we will assume the following.
\begin{description}
\item
\begin{condition}
\label{R0}  $\rho$ and $\rho_{0}$ are bounded $\rho$-functions.
\end{condition}
\end{description}

Define the residual scale functional $S^{\ast}(G,\xi)$ by%
\begin{equation}
\mathrm{E}_{G}\rho_{0} \biggl(  \frac{y-\underline{g}(x,\xi)}{S^{\ast}(G,\xi
)} \biggr)  =\delta, \label{defS*}%
\end{equation}
with $\delta\in(0,1)$. Then the regression S functional $T_{\mathrm{S}}$ and
the associated error scale M functional $S(G)$ are, respectively, defined by
\begin{equation}
T_{\mathrm{S}}(G):= (  T_{\mathrm{S},\beta}(G),T_{\mathrm{S},\alpha
}(G) )  =\arg\min_{\xi\in B\times R}S^{\ast}(G,\xi) \label{S}%
\end{equation}
and%
\begin{equation}
S(G)=\min_{\xi\in B\times R}S^{\ast}(G,\xi). \label{defSG}%
\end{equation}

We will deal with a regression M functional $T_{\mathrm{M}}(G)$ defined as
\begin{equation}
T_{\mathrm{M}}(G):= (  T_{\mathrm{M},\beta}(G),T_{\mathrm{M},\alpha
}(G) )  =\arg\min_{\xi\in B\times R}M_{G}(\xi), \label{M-loc}%
\end{equation}
where the function $M_{G}\dvtx B\times R\rightarrow R$ is%
\begin{equation}
M_{G}(\xi)=\mathrm{E}_{G}\rho \biggl(  \frac{y-\underline{g}(x,\xi)}%
{\widetilde{S}(G)} \biggr)  \label{MG}%
\end{equation}
and $\widetilde{S}(G)$ is an arbitrary residual scale functional, for
example, the one defined in (\ref{defSG}).

It is easy to show that the S regression functional defined in (\ref{S}) is
also an M functional. In fact, $T_{\mathrm{S}}(G)$ coincides with
$T_{\mathrm{M}}(G)$ when  in  (\ref{MG}) we have $\rho=\rho_{0}$ and
$\widetilde{S}(G)=S(G).$ We may then write%
\begin{equation}
T_{\mathrm{S}}(G)=\arg\min_{\xi\in B\times R}\mathrm{E}_{G}\rho_{0} \biggl(
\frac{y-\underline{g}(x,\xi)}{S(G)} \biggr)  . \label{s=m}%
\end{equation}

\begin{remark}
\label{RemaUni}
In general, the minimum at (\ref{S}) or (\ref{M-loc}) might be
attained at more than one value of $\xi$. It will be henceforth assumed that
the functional is well-defined by the choice of a single value. Our results
will not depend on how the choice is made. However, it will be shown in
Section~\ref{SecFish-M} that  under very general conditions,  if $G_{0}$
is  the distribution of $(x,y)$ satisfying (\ref{modelogral}), then
$T_{\mathrm{S}}(G_{0})$ and $T_{\mathrm{M}}(G_{0})$ are unique  and
$T_{\mathrm{S,}\beta}(G_{0}) =T_{\mathrm{M,}\beta}(G_{0})=\beta_{0}$ (Fisher-consistency).
\end{remark}

\section{Fisher-consistency of M and S estimates}\label{SecFish-M}

In this section, we give sufficient conditions to guarantee that both (\ref{S})
and (\ref{M-loc}) are minimized at unique values, and that $T_{\mathrm{M},\beta
} (G_{0})=T_{\mathrm{S}, \beta} (G_{0})=\beta_{0}$ (Fisher consistency for~$\beta_{0}$).

Recall that a density $f$ is \emph{strongly unimodal}   if there exists $a$
such that $f (  t )  $ is nondecreasing for $t<a,$ nonincreasing for
$t>a,$ and $f$ has a unique maximum at $t=a$.

We will need the following condition on $\rho$.

\begin{condition}
\label{R1} The function $\rho$ is a $\rho$-function such that  for some
$m>0,$ $\rho(u)=1$ iff $|u|\geq m,$ and $\log(1-\rho)$ is concave on
$(-m,m)$.
\end{condition}

It is easy to check that Condition~\ref{R1} with $m=k$ holds in particular for
the popular family of bisquare functions, defined by%
\[
\rho_{k} (  u )  =1- \biggl(  1- \biggl(  \frac{u}{k} \biggr)
^{2} \biggr)  ^{3}I(|u|\leq k).
\]

We will establish the Fisher-consistency of $T_{\mathrm{M}}$. Put for brevity
$\sigma=S (  G_{0} )  $ and let $F_{0}$ be the distribution of
$u_{i}$ in (\ref{modelogral}) and assume that it has a strongly unimodal
density. Let $\Delta$ denote the unique minimizer of $\mathrm{E}_{F_{0}}%
\rho( (  u-t )  /\sigma)$; note that if $u_{i}$ is symmetric around
$\mu_{0}$, then part (b) of Theorem~\ref{Teore4} implies that $\Delta=\mu
_{0}.$

\begin{theorem}
\label{Teore5}  Let $G_{0}$ be the joint distribution of  $(x_{i},y_{i})$
satisfying model  (\ref{modelogral}), where $u_{i}$ has distribution
$ F_{0}$ with a strongly unimodal density. Assume that Conditions
\ref{CondId} and~\ref{R1} hold. Then $M_{G_{0}}(\xi)$ is minimized at the
unique point $T_{\mathrm{M}} (  G_{0} )  = (  \beta_{0}%
,\Delta )  $, and so $T_{\mathrm{M}}$ is Fisher-consistent for
$\beta_{0}$, that is, $T_{\mathrm{M},\beta}(G_{0})=\beta_{0}$. If we also assume that
$F_{0}$ is symmetric around $\mu_{0},$ we have $T_{\mathrm{M},\alpha}(G_{0})=\mu_{0}.$
\end{theorem}

\begin{remark}
\label{Rt5} Theorem~\ref{Teore5} gives also sufficient conditions for the
Fisher-consistency of the regression S functional $T_{\mathrm{S}}$. In
fact$,$ according to (\ref{s=m}), $T_{\mathrm{S}}$ is also an M functional.
\end{remark}

\section{Weak continuity of M and S regression functionals}\label{SecConti-M-S}

\begin{definition}
We say that a functional $T$ is \emph{weakly continuous} at $G$ if
$G_{n}\rightarrow_{w}G$ implies $T(G_{n})\rightarrow T(G).$
\end{definition}

We will show the weak continuity of the functionals defined above in two
cases: nonlinear regression with a compact parameter space $B$, and linear regression.

Define for $G=\mathcal{D} (  x,y )  $
\begin{equation}
c(G)=\sup\{\mathrm{P}_{G}(\beta^{\prime}x+\alpha=0)\dvt\beta\in R^{p},\alpha\in
R,\beta\neq0\}. \label{c}%
\end{equation}

\begin{theorem}
\label{Teore01}Let $G_{0}=\mathcal{D} (  x,y )  $ be such that
(\ref{M-loc}) has a unique solution $T_{\mathrm{M}}(G_{0})$.  Assume that
$\widetilde{S}$ is weakly continuous at $G_{0}$ and $\widetilde{S}(G_{0})>0$.
Then $T_{\mathrm{M}}=(T_{\mathrm{M},\beta},T_{\mathrm{M},\alpha})$  is
weakly continuous at $G_{0}$ if  either  \textup{(a)} or \textup{(b)} holds, where
\begin{enumerate}[(b)]
\item[(a)] $B$  is compact,
\item[(b)] $B=R^{p}$,  $g(x,\beta)=\beta^{\prime}x$ and
\begin{equation}
M_{G_{0}}(T_{\mathrm{M}}(G_{0}))<1-c(G_{0}). \label{M<1}%
\end{equation}
\end{enumerate}
\end{theorem}

\begin{theorem}
\label{Teore02}Let $G_{0}=\mathcal{D} (  x,y )  $ be such that
$T_{\mathrm{S}}(G_{0})$ is unique and $S(G_{0})>0.$ Assume that either \textup{(a)}
$B$ is compact, or \textup{(b)} $B=R^{p}$, $g$ is linear, that is, $g(x,\beta
)=\beta^{\prime}x$  and $\delta<1-c(G_{0})$ with $c(G)$ defined in
(\ref{c}). Then $S(G)$ and $T_{\mathrm{S}}(G)=(T_{\mathrm{S},\beta
},T_{\mathrm{S},\alpha})$ are weakly continuous at $G_{0}$.
\end{theorem}

Let now $G_{0}$ be the distribution of $ (  x,y )  $ under model
(\ref{modelogral}), and assume that $T_{\mathrm{M}}$ (resp.,
$T_{\mathrm{S}}$) is Fisher-consistent for $\beta_{0},$ that is, $T_{\mathrm{M}%
,\beta} (  G_{0} )  =\beta_{0}$ (resp., $T_{\mathrm{S},\beta
} (  G_{0} )  =\beta_{0}).$ Then the former results imply that
$T_{\mathrm{M},\beta}$ (resp., $T_{\mathrm{S},\beta})$ evaluated at the
empirical distribution is consistent whenever the empirical distributions
converge to the underlying one. More precisely, we have the following result.

\begin{corollary}
\label{CoroConsiste} Assume the same hypotheses as in Theorem~\ref{Teore01}
(resp., Theorem~\ref{Teore02}) plus the Fisher-consistency of
$T_{\mathrm{M}}$ (resp., $T_{\mathrm{\mathrm{S}}}$): $T_{\mathrm{M}%
,\beta}(G_{0})=T_{\mathrm{S},\beta}(G_{0})=\beta_{0}$. Call $G_{n}$ the
empirical distribution of $\{ (  x_{i},y_{i} )  \dvt i=1,\ldots,n\}$. If
$G_{n}\rightarrow_{w}G_{0}$ a.s., then $\{T_{\mathrm{M},\beta}(G_{n})\}$
(resp., $\{T_{\mathrm{S},\beta}(G_{n})\}$) is strongly consistent for
$\beta_{0}$.
\end{corollary}

This result is immediate. The a.s. weak convergence of $G_{n}$ to $G_{0}$ is
well known to hold for i.i.d. $ (  x_{i},y_{i} )  $. It holds also
under more general assumptions on the joint distribution of $\{ (
x_{i},y_{i} )  \dvt i\geq1\},$ such as ergodicity.

We now turn to qualitative robustness. Consider a sequence of estimates
$\{\widehat{\xi}_{n}\}$ based on a functional $T,$ that is, $\widehat{\xi}%
_{n}=T (  G_{n} )  $ where $G_{n}$ is the empirical distribution
corresponding to data $ (  z_{1},\ldots,z_{n} )  .$ Hampel~\cite{Hampel}
proved that for $\{\widehat{\xi}_{n}\}$ to be qualitatively robust at a
distribution $G_{0}$ it suffices that $T$ be weakly continuous at $G_{0}$ and
$\widehat{\xi}_{n}$ be a continuous function of $ (  z_{1},\ldots,z_{n}%
 )  .$

Papantoni-Kazakos and Grey~\cite{Papantoni} employ a weaker definition of
robustness, which they call \emph{asymptotic qualitative robustness}, and
prove that it is equivalent to weak continuity. Therefore, Theorems
\ref{Teore01} and~\ref{Teore02} imply the asymptotic qualitative robustness of
$T_{\mathrm{M}}$ and $T_{\mathrm{S}}.$

\section{Differentiability of estimating functionals}\label{SecDife}

 In this section, we shall first deal with the differentiability of general
functionals and then specialize to our regression case. Let $\mathcal{G}_{h}$
 be a set of  distributions on $R^{h}.$ Consider an estimating functional
$T\dvtx \mathcal{G}_{h} \rightarrow R^{k}$. Hampel~\cite{Hampel2} defines the
\emph{influence function} of $T$   at $G\in\mathcal{G}_{h}$ as the function
$I_{T,G}(z )  \dvtx R^{h}\rightarrow R^{k}$
\begin{equation}
I_{T,G}(z )=    \frac{\partial(T((1-\varepsilon)G+\varepsilon
\delta_{z}))}{\partial\varepsilon} \bigg| _{\varepsilon=0}, \label{inffund}%
\end{equation}
where $\delta_{z}$ is the point mass distribution at $z$. Given a distance $d$
on $\mathcal{G}_{h}$ which metricizes the topology of convergence in
distribution, $T$ is \emph{Fr\'{e}chet differentiable}   at $G_{0}$ under $d$
if%
\[
T(G)-T(G_{0})=\mathrm{E}_{G}I_{T,G_{0}}(z ) +\mathrm{o}(d(G,G_{0})).
\]

Fr\'{e}chet differentiability can be used to prove the asymptotic normality of
the estimate. However, Fr\'{e}chet differentiability also requires that
$I_{T,G}(z ) $ be bounded. Since this condition  is not satisfied by
regression M estimates, we are going to define a weaker type of
differentiability, which suffices to prove asymptotic normality.

\begin{definition}
\label{DefiDife}
Let $T$ be  an estimating functional that is weakly
continuous at $G_{0}.$ We say that $T$ is \emph{weakly differentiable} at
 a sequence $\{G_{n}\}$ converging weakly to $G_{0}$ if
\begin{equation}
T(G_{n})-T(G_{0})=\mathrm{E}_{G_{n}}I_{T,G_{0}}(z  )+ \mathrm{o}(
 \Vert \mathrm{E}_{G_{n}}I_{T,G_{0}}(z {)} \Vert  )  .
\label{stdif}%
\end{equation}
\end{definition}

The definition of weak differentiability helps understanding the asymptotic
behavior of $T(G_{n})-T(G_{0})$, as the next lemma shows.

\begin{lemma}
\label{lem1} Consider a random sequence of distributions $\{G_{n}\}$
converging weakly to $G_{0}$ a.s. Assume that $T$ is weakly differentiable at
$\{G_{n}\}$ a.s. and that for some sequence $\{a_{n}\}$ of real numbers
\[
a_{n}\mathrm{E}_{G_{n}}I_{T,G_{0}}(z )\rightarrow _{d}H.
\]
Then
\begin{equation}
a_{n}\bigl(T(G_{n})-T(G_{0})\bigr)=a_{n}\mathrm{E}_{G_{n}}I_{T,G_{0}}(z )+
\mathrm{o}_{p}(1) \label{exp1}%
\end{equation}
and therefore $a_{n}(T(G_{n})-T(G_{0}))\rightarrow_{d}H$ too.
\end{lemma}

The proof of this lemma is immediate.

\begin{remark*}
Note that if (\ref{exp1}) holds for a joint functional
$T=(T_{1},T_{2}),$ it also holds for $T_{1},$ that is,
\begin{equation}
a_{n}\bigl(T_{1}(G_{n})-T_{1}(G_{0})\bigr)=a_{n}\mathrm{E}_{G_{n}}I_{T_{1},G_{0}%
}(z )+ \mathrm{o}_{p}(1).
\end{equation}
\end{remark*}

We now deal with the differentiability of a \emph{general M estimating
functional}, that is, a functional $T$ defined on a subset of $\mathcal{G}_{p}$
with values in $ R^{q}$, that for some function $\Psi\dvtx R^{p}\times
R^{q}\rightarrow R^{q}$ satisfies the equation%
\begin{equation}
\mathrm{E}_{G}\Psi(z,T(G))=0. \label{Mfuncgen}%
\end{equation}

We will assume that $\Psi$ is continuously differentiable with respect to
$\theta$ and call $\dot{\Psi}(z,\theta)$ (or alternatively $\partial{\Psi
}(z,\theta)/\partial\theta$) the $q\times q$ differential  matrix with
elements $\dot{\Psi}_{jk}(z,\theta)=\partial\Psi_{j}(z,
\theta)/\partial\theta_{k}$. Define%
\begin{equation}
D (  G,\theta )  =\mathrm{E}_{G}\dot{\Psi} (  z,\theta )  .
\label{defD}%
\end{equation}

Let $\theta_{0}=T(G_{0})$ and assume that%
\begin{equation}
D_{0}=D (  G_{0},\theta_{0} )  \label{defD0}%
\end{equation}
exists. Assume that $T$ is weakly continuous at $G_{0}$ and that the following holds.

\begin{condition}
\label{CondPsi}
$D_{0}$ is nonsingular and there exists $\eta>0$ such that
\begin{equation}
\mathrm{E}_{G_{0}}\sup_{\|\theta-\theta_{0}\|\leq\eta}\|\dot{\Psi} (
z,\theta )  \|<\infty, \label{SupPsiDot}%
\end{equation}
 where $ \Vert \cdot \Vert $ denotes the $l_{2}$ norm.
\end{condition}

Then, it is easy to show that the influence function of $T$ at $G_{0}$ is
given by%
\begin{equation}
I_{T,G_{0}}(z )=- D_{0}^{-1}\Psi(z,\theta_{0}). \label{IF-Mest}%
\end{equation}

We shall now see that the following conditions are sufficient for the weak
differentiability of $T$ at $\{G_{n}\}$.

\begin{condition}
\label{CondiA} $\{G_{n}\}$ is a sequence of distribution functions that
converges weakly to $G_{0}$ and
\begin{equation}
\lim_{\eta\rightarrow0}\limsup_{n\rightarrow\infty}\sup_{ \Vert
\theta-\theta_{0} \Vert \leq\eta} \Vert D(G_{n},\theta)-D_{0}%
 \Vert =0. \label{cca}%
\end{equation}
\end{condition}

\begin{condition}
\label{CondiC}
$\{G_{n}\}$ is a sequence of distribution functions such that,
at a neighborhood of $\theta_{0}$, for each $n$
\begin{equation}
\frac{\partial}{\partial\theta}\mathrm{E}_{G_{n}}\Psi(z,%
\theta)=\mathrm{E}_{G_{n}}\frac{\partial}{\partial\theta}\Psi(z,%
\theta). \label{conmuta}%
\end{equation}
\end{condition}

Condition~\ref{CondiA} means that $D(G_{n},\theta)$ approaches $D(G_{0}%
,\theta_{0})$ when $n$ is large and $\theta$ is close to~$\theta_{0}.$
Condition~\ref{CondiC} means that we can interchange differentiation of
$\Psi(z,\theta)$ with respect to $\theta$ and expectation  with
respect to $G_{n}$. Theorem~\ref{TeoDiffe} shows that these two conditions
imply weak differentiability and Theorem~\ref{Teorevarios} shows that these
conditions hold in very general circumstances.

\begin{theorem}
\label{TeoDiffe}
Assume that $ T$  is an M functional satisfying
(\ref{Mfuncgen}) and weakly continuous at $G_{0},$ that $\dot{\Psi
}(z,\theta)$ is continuous in $\theta,$ and that Condition
\ref{CondPsi} holds. If $\{G_{n}\}$ satisfies Conditions~\ref{CondiA} and
\ref{CondiC}; then  $T$ is weakly differentiable at $\{G_{n}\}.$
\end{theorem}

The following theorem gives sufficient conditions for a.s. differentiability
 of M functionals, at a random sequence of distributions.

\begin{theorem}
\label{Teorevarios} Let $\{G_{n}\}$ be a sequence of random distributions
converging weakly to $G_{0}$ and satisfying Condition~\ref{CondiC}
a.s. Assume also that $\dot{\Psi} (  z,\theta )  $ is continuous in
$\theta$ and that Condition~\ref{CondPsi} holds. Let $ T$ be an M functional
satisfying (\ref{Mfuncgen}) and weakly continuous at $G_{0}$. Then $T$ is
weakly differentiable at $\{G_{n}\}$ a.s. in any of the following two cases:
\textup{(a)}~for each function $d(z)$ such that $\mathrm{E}_{G_{0}} \vert
d(z) \vert <\infty$, on a set of probability one we have that
$\{\mathrm{E}_{G_{n}}d(z)\}$ converges to $\mathrm{E}_{G_{0}}d(z)$, or \textup{(b)}~$\dot{\Psi}(z,\theta)$ is bounded.
\end{theorem}

Note case (a) contains situations where a Law of Large Numbers holds, in
particular when $G_{n}$ is the empirical distribution of an ergodic process.

\begin{corollary}
Let $\{G_{n}\}$  be a sequence of empirical distributions associated to
i.i.d. $\{z_{i}\}$ with distribution $G_{0}.$ Assume that $\dot{\Psi} (
z,\theta )  $ is continuous in $\theta$, that Condition~\ref{CondPsi}
holds and that $I_{T,G_{0}} (  z )  $ has finite second moments under
$G_{0}$. Let $T$ be an M functional continuous at $G_{0}$. Then $n^{1/2}%
 (  T (  G_{n} )  -T (  G_{0} )   )  \rightarrow
_{d}N (  0,V )  $ with
\begin{equation}
V=\mathrm{E}_{G_{0}}I_{T,G_{0}} (  z )  I_{T,G_{0}} (  z )
^{\prime}. \label{defV}%
\end{equation}
\end{corollary}

There are many examples where Fr\'{e}chet differentiability does not hold and
that can  be dealt with using the concept of weak differentiability.  One of
these cases is that of MM estimates for linear an nonlinear regression which
is treated in detail in the next section. Other examples where Fr\'{e}chet
differentiability fails are MM estimates for the multivariate linear model
(see~\cite{Kudraszov}) and M estimates for logistic
models (see~\cite{Bianco} and~\cite{Croux}). An example where the asymptotic expansions  that can be
obtained with weak differentiability are essential to prove asymptotic
normality is the problem of robust estimation with missing data considered
by Sued and Yohai~\cite{Sued-Yohai}.

\section{MM estimates}\label{SectionMM}

In this section, we will summarize the properties derived from Theorems 1--6
for S and MM estimates of regression and location.

\subsection{Regression case}

Recall that MM estimates, which we denote here by $T_{\mathrm{MM}%
}=(T_{\mathrm{MM},\beta},T_{\mathrm{MM},\alpha}),$ are defined in
(\ref{M-loc}), where $\widetilde{S}$ is the functional $S$  defined in
(\ref{defSG}). For notational convenience, we shall call $\rho_{1}$ the
$\rho$-function employed in (\ref{MG}), and we will assume that $\rho_{1}%
\leq\rho_{0}$. As mentioned above, the definition of $\widehat{\xi
}_{\mathrm{MM}}$ in (\ref{defMestima}) requires also $\widehat{\sigma}$
defined by (\ref{Sscale}), and hence also $\widehat{\xi}_{\mathrm{S}}$ defined
in (\ref{defSestima}). Therefore,  these three estimates must be considered
simultaneously. Call
\begin{equation}
\widehat{\theta}= (  \widehat{\xi}_{\mathrm{S}},\widehat{\xi}%
_{\mathrm{MM}},\widehat{\sigma} )  \label{joint}%
\end{equation}
the joint solution of (\ref{defMestima})--(\ref{defSestima})--(\ref{Sscale}).

In the rest of this  section, we assume the following properties.

\begin{condition}
\label{R2} $\rho_{0}$   and $\rho_{1}$ are twice continuously
differentiable.
\end{condition}

We denote by $\psi_{0}$ and $\psi_{1}$ the derivatives of $\rho_{0}$ and
$\rho_{1},$ respectively. Assume also the following condition.

\begin{condition}
\label{R3} $g$ is twice continuously differentiable with respect to
$\beta.$
\end{condition}

We denote  by $\underline{\dot{g}}(x,\xi)$ and $\underline{\ddot{g}} (
x,\xi )  $ the vector of first derivatives and the matrix of second
derivatives of $\underline{g}$ with respect to $\xi$, respectively.
Analogously, we denote by $\dot{g}(x,\beta)$ and $\ddot{g} (
x,\beta )  $ the vector of first derivatives and the matrix of second
derivatives of $g$ with respect to $\beta$, respectively. Note that
$\underline{\dot{g}}(x,\xi)$ and $\underline{\ddot{g}} (  x,\xi )  $
depend only on $\beta$, and for this reason we will indistinctly use also the
notation $\underline{\dot{g}}(x,\beta)$ and $\underline{\ddot{g}} (
x,\beta )  .$

Differentiating (\ref{defMestima}) we have that $\widehat{\xi}_{\mathrm{MM}%
}$ satisfies the system%
\begin{equation}
\frac{1}{n}\sum_{i=1}^{n}\psi_{1} \biggl(  \frac{y_{i}-\underline{g}%
(x_{i},\widehat{\xi}_{\mathrm{MM}})}{\widehat{\sigma}} \biggr)  \underline
{\dot{g}}(x_{i},\widehat{\xi}_{\mathrm{MM}})=0. \label{sys1}%
\end{equation}
It is immediate that $\widehat{\xi}_{\mathrm{S}}$ also satisfies
\[
\widehat{\xi}_{\mathrm{S}}=\arg\min_{\xi\in B\times R}\frac{1}{n}\sum
_{i=1}^{n}\rho_{0} \biggl(  \frac{y_{i}-\underline{g}(x_{i},\xi)}{\widehat
{\sigma}} \biggr)  .
\]
Then, differentiating this equation we get
\begin{equation}
\frac{1}{n}\sum_{i=1}^{n}\psi_{0} \biggl(  \frac{y_{i}-\underline{g}%
(x_{i},\widehat{\xi}_{\mathrm{S}})}{\widehat{\sigma}} \biggr)  \underline
{\dot{g}}(x_{i},\widehat{\xi}_{\mathrm{S}})=0. \label{sys2}%
\end{equation}

Finally according to (\ref{defMscale}), $\widehat{\sigma}$ satisfies%
\begin{equation}
\frac{1}{n}\sum_{i=1}^{n}\rho_{0} \biggl(  \frac{y_{i}-\underline{g}%
(x_{i},\widehat{\xi}_{\mathrm{S}})}{\widehat{\sigma}} \biggr)  -\delta=0.
\label{sys3}%
\end{equation}
Then $\widehat{\theta}$ satisfies the system of $2q+3$ equations
(\ref{sys1})--(\ref{sys3}). Putting $z_{i}=(x_{i},y_{i})$ and
denoting by $G_{n}$ the empirical distribution of $\{z_{1},\ldots,z_{n}\},$ this
system can be written as%
\begin{equation}
\frac{1}{n}\sum_{i=1}^{n}\Psi (  z_{i},\widehat{\theta} )
=\mathrm{E}_{G_{n}}\Psi (  z,\widehat{\theta} )  =0,
\label{EstimaEquat}%
\end{equation}
where if $\theta= (  \xi_{\mathrm{S}},\xi_{\mathrm{MM}},\sigma )  ,$
$\Psi(z,\theta)$ is defined by%
\[
\everymath{\displaystyle}
\Psi (  z,\theta )  = \left[
\begin{array}{@{}c@{}}%
\psi_{0} \biggl(  \frac{y-\underline{g}(x,\xi_{\mathrm{S}})}{\sigma} \biggr)
\underline{\dot{g}}(x,\xi_{\mathrm{S}})\\[10pt]
\psi_{1} \biggl(  \frac{y-\underline{g}(x,\xi_{\mathrm{MM}})}{\sigma} \biggr)
\underline{\dot{g}}(x,\xi_{\mathrm{MM}})\\[10pt]
\rho_{0} \biggl(  \frac{y-\underline{g}(x,\xi_{\mathrm{S}})}{\sigma} \biggr)
-\delta.
\end{array}
 \right]  .
\]

Let%
\begin{equation}
T(G)= (  T_{\mathrm{S}} (  G )  ,T_{\mathrm{MM}} (  G )
,S (  G )   )  \label{TTT}%
\end{equation}
be the  estimating functional  associated to $\widehat{\theta}.$  Then, if
(\ref{SupPsiDot}) holds, we can differentiate the functions to be minimized in
 (\ref{M-loc}) and (\ref{s=m}) inside the expectation, obtaining that
$ T(G)$  satisfies the equation%
\begin{equation}
\mathrm{E}_{G}\Psi(z,T(G))=0. \label{Tcompleto}%
\end{equation}
Note that the solution to this equation is in general not unique, and
therefore, $T$ is not defined exclusively by the equation.

To verify (\ref{SupPsiDot}), in addition to Conditions~\ref{R0}--\ref{R1}%
--\ref{R2}--\ref{R3} we need the following assumption:

\begin{condition}
\label{R4} For some $\eta>0 $%
\begin{equation}
\mathrm{E}_{G_{0}}\sup_{ \Vert \beta-\beta_{0} \Vert \leq\eta
} \Vert {\dot{g}}(x,\beta) \Vert ^{2}<\infty  \quad \mbox{and} \quad
  \mathrm{E}_{G_{0}} \sup_{ \Vert \beta-\beta_{0} \Vert \leq\eta
} \Vert {\ddot{g}}(x,\beta) \Vert <\infty. \label{conDife-g}%
\end{equation}
Assume that $D_{0}$ defined by (\ref{defD0}) is nonsingular; then under these
assumptions, we also get that $I_{T,G_{0}} (  z )  $ has finite
second moments under $G_{0}.$ Note that in the case of linear regression,
(\ref{conDife-g}) reduces to $\mathrm{E}_{G_{0}} \Vert x \Vert
^{2}<\infty.$
\end{condition}

Define%
\begin{equation}
\alpha_{0i}=\arg\min_{t}\mathrm{E}_{F_{0}}\rho_{i} \biggl(  \frac{u-t}{S(G_{0}%
)} \biggr)  , \qquad  i=0,1, \label{delta01}%
\end{equation}
where $F_{0}$ is the distribution of $u_{i}$ in model (\ref{modelogral}). We
will see in Theorem~\ref{teoremaso} that under some general conditions,
$T_{\mathrm{S,}\alpha}(G_{0})=\alpha_{00}$ and $T_{\mathrm{MM},\alpha}%
(G_{0})=\alpha_{01.}.$

Put $\theta_{0}=(\beta_{0},\alpha_{00},\beta_{0},\alpha_{01},\sigma_{0})$ with
$\sigma_{0}=S(G_{0}).$ The following numbers, vectors and matrices are required
to derive a closed formula for the influence functions of $T_{{\mathrm{MM}}}$
and~$T_{\mathrm{S}}$. Let
\begin{eqnarray*}
a_{0i}&=&\mathrm{E}_{G_{0}}\psi_{i}^{\prime} \biggl(  \frac{y-g(x,\beta
_{0})-\alpha_{0i}}{\sigma_{0}} \biggr)  =\mathrm{E}_{F_{0}}\psi_{i}^{\prime
} \biggl(  \frac{u-\alpha_{0i}}{\sigma_{0}} \biggr)  , \qquad  i=0,1,
\\
e_{0i}&=&\mathrm{E}_{F_{0}}\biggl (  \frac{u-\alpha_{0i}}{\sigma_{0}} \biggr)
\psi_{0}^{\prime} \biggl(  \frac{u-\alpha_{0i}}{\sigma_{0}} \biggr)  , \qquad  i=0,1,
\\
d_{0}&=&\mathrm{E}_{F_{0}} \biggl(  \frac{u-\alpha_{00}}{\sigma_{0}} \biggr)
\psi_{0} \biggl(  \frac{u-\alpha_{00}}{\sigma_{0}} \biggr)  ,
\\
b_{0}&=&\mathrm{E}_{G_{0}}\dot{g}(x,\beta_{0}), \qquad   b_{0}^{\ast}=(b_{0}^{\prime
},1)^{\prime},
\\
A_{0}&=&\mathrm{E}_{F_{0}}\bigl(\dot{g}(x,\beta_{0})-b_{0}\bigr)\bigl(\dot{g}(x,\beta
_{0})-b_{0}\bigr)^{\prime}%
\end{eqnarray*}
and%
\begin{equation}\label{defC0}
C_{0}= \left[
\begin{array}
[c]{@{}c@{ \quad }c@{}}%
A_{0}+b_{0}b_{0}^{\prime} & b_{0}\\
b_{0}^{\prime} & 1
\end{array}
 \right]  .
\end{equation}

It is shown in Section~\ref{secDerIF} that the influence function of
$T_{\mathrm{MM}}$ is given by%
\begin{equation}
I_{T_{\mathrm{MM},\beta},G_{0}}(x,y)=\frac{\sigma_{0}}{a_{01}}\psi_{1} \biggl(
\frac{y-\underline{g}(x,(\beta_{0},\alpha_{01}))}{\sigma_{0}} \biggr)
A_{0}^{-1} \bigl(  \dot{g}(x,\beta_{0})-b_{0} \bigr)  \label{INFMM1}%
\end{equation}
and%
\begin{eqnarray}\label{INFMM2}
I_{T_{\mathrm{MM},\alpha},G_{0}}(x,y)  &  =&-\frac{\sigma_{0}}{a_{01}}\psi
_{1} \biggl(  \frac{y-\underline{g}(x, (  \beta_{0},\alpha_{01} )
)}{\sigma_{0}} \biggr)\bigl   [  1+b_{0}^{\prime}A_{0}^{-1} \bigl(  b_{0}-\dot
{g}(x,\beta_{0}) \bigr)   \bigr]\nonumber
\\[-8pt]
\\[-8pt]
&&{}  +\frac{\sigma_{0}e_{01}}{a_{01}d_{0}} \biggl(  \rho_{0} \biggl(  \frac
{y-\underline{g}(x,(\beta_{0},\alpha_{01}))}{\sigma_{0}} \biggr)
-\delta \biggr)  . %
\nonumber
\end{eqnarray}

The influence functions of $ T_{\mathrm{S},\beta}$ and $T_{\mathrm{S}%
,\alpha}$ can be obtained similarly replacing  $\alpha_{01}, a_{01}$ and
$e_{01}$ by $\alpha_{00}, a_{00}$ and $e_{00}$, respectively.

If the errors $u_{i}$ have a symmetric distribution $F_{0}$, then $e_{01}=0$
and $\alpha_{01}=\alpha_{00}=\alpha_{0},$ the center of symmetry of $F_{0}%
.$ This entails a considerable simplification of the influence function
 $I_{T_{\mathrm{MM}}}$. In fact, in this case, we get%
\begin{equation}
I_{T_{\mathrm{MM}},G_{0}}(z )= \frac{\sigma_{0}}{\mathrm{E}_{F_{0}%
}\psi_{1}^{\prime} (  (u-\alpha_{0})/\sigma_{0} )  }\psi_{1} \biggl(
\frac{y-g(x,\beta_{0})-\alpha_{0}}{\sigma_{0}} \biggr)  C_{0}^{-1}%
\underline{\dot{g}}(x,\beta_{0}), \label{INFMM4}%
\end{equation}
 and the asymptotic covariance matrix (\ref{defV}) is%
\begin{equation}
V=\sigma_{0}^{2}\frac{\mathrm{E}_{F_{0}}\psi_{1} (  (u-\alpha_{0}%
)/\sigma_{0} )  ^{2}}{ (  \mathrm{E}_{F_{0}}\psi_{1}^{\prime} (
(u-\alpha_{0})/\sigma_{0} )   )  ^{2}}C_{0}^{-1}. \label{varsym}%
\end{equation}

The next theorem summarizes the properties of S and MM  regression functionals.\vspace*{-2pt}

\begin{theorem}
\label{teoremaso} Let $z=(x,y)$ satisfy model (\ref{modelogral}) where the
distribution $F_{0}$ of $u_{i}$ has a strongly unimodal density and Condition
\ref{CondId} holds.  Assume that $\rho_{0}$ and $\rho_{1}$ are bounded $\rho
$-functions that satisfy Condition~\ref{R1}, with $\rho_{1}(u)\leq\rho_{0}%
(u)$. Let $T$ be defined by (\ref{TTT}) and let $G_{0}$ be the
distribution of $(x,y)$. Then:
\begin{enumerate}[(viii)]
\item[(i)] $T_{\mathrm{S}}(G_{0})=(\beta_{0},\alpha_{00})$ is the unique
minimizer in (\ref{S}). If $F_{0}$ is symmetric with respect to $\mu_{0}$, then
$\alpha_{00}=\mu_{0.}$

\item[(ii)] $T_{\mathrm{MM}}(G_{0})=(\beta_{0},\alpha_{01})$ is the unique
minimizer in (\ref{M-loc}).  If $F_{0}$ is symmetric with respect to $\mu
_{0}$ then $\alpha_{01}=\mu_{0}$.

\item[(iii)] The functional $T=(T_{\mathrm{S}},T_{\mathrm{MM}},S)$ is weakly
continuous at $G_{0}$ if either \textup{(a)} $B$ is compact, or \textup{(b)} $B=R^{p},$
$g(x,\beta)=\beta^{\prime}x$ and $\delta<1-c(G_{0})$.

\item[(iv)] Assume also that Conditions~\ref{R2},~\ref{R3} and~\ref{R4}
 hold, that $a_{00}\neq0,$ $a_{01}\neq0,$ $d_{0}\neq0$ and that $A_{0}$ is
invertible. Then, $D_{0}=\mathrm{E}_{G_{0}}\dot{\Psi} (  z,T(G_{0}%
) )  $ is invertible, $I_{T_{\mathrm{MM,}\beta},G_{0}}(x,y)$ and
$I_{T_{\mathrm{MM,}\alpha},G_{0}}(x,y)$ are given by (\ref{INFMM1}) and
(\ref{INFMM2}), respectively, while the influence functions $I_{T_{\mathrm{S,}%
\beta},G_{0}}(x,y)$ and $I_{T_{S,\alpha},G_{0}}(x,y)$ have a similar
expression replacing  $\alpha_{01}, a_{01}$ and $e_{01}$ by $\alpha
_{00}, a_{00}$ and $e_{00}$, respectively.

\item[(v)] Under the same assumptions as in \textup{(iv)}, let $\{G_{n}\}$ be a
sequence of random distributions converging weakly to $G_{0}$ and satisfying
Condition~\ref{CondiC} a.s. Suppose also that for each function $d(z)$ such
that $\mathrm{E}_{G_{0}} \vert d(z) \vert <\infty$,  we have that
$\{\mathrm{E}_{G_{n}}d(z)\}$ converges to $\mathrm{E}_{G_{0}}d(z)$ a.s. Then,
the functional $T$ is weakly differentiable at $\{G_{n}\}$.

\item[(vi)] Assume the same conditions as in \textup{(v)} and:%
\begin{equation}
n^{1/2}\mathrm{E}_{G_{n}}I_{T,G_{0}}(x,y)\rightarrow_{d}H. \label{convH}\vspace*{-2pt}
\end{equation}
Then\vspace*{-2pt}
\begin{equation}
n^{1/2}\bigl(T(G_{n})-T(G_{0})\bigr)=n^{1/2}\mathrm{E}_{G_{n}}I_{T,G_{0}}(x,y)+\mathrm{o}_{p}(1)
\label{asinMM}\vspace*{-2pt}
\end{equation}
and therefore
\begin{equation}
n^{1/2}\bigl(T(G_{n})-T(G_{0})\bigr)\rightarrow_{d}H. \label{iv1}%
\end{equation}

\item[(vii)] Assume that the conditions in \textup{(iv)} hold and that $\{(x_{i}%
,u_{i})\dvt i\geq1\}$  are i.i.d. Let $ G_{n}$ be the sequence of empirical
distributions corresponding to  \mbox{$\{(x_{i},y_{i})\dvt i\geq1\}$} with common
distribution $G_{0}$. Then (\ref{iv1}) holds with $H =N(0,V)$ and
$V=\break\mathrm{E} [  I_{T_{\mathrm{MM}},G_{0}}(x,y)I_{T_{\mathrm{MM}},G_{0}%
}(x, y)^{\prime} ]  ,$ where $I_{T_{\mathrm{MM}},G_{0}}(x,y)$ is
defined by (\ref{INFMM1}) and~(\ref{INFMM2}).

\item[(viii)] Assume that the conditions in \textup{(iv)} hold, that $\{u_{i}\dvt i\geq1\}$
is stationary and ergodic and that $\{x_{i}, i\geq1\}$ are i.i.d. and
independent of $\{u_{i}\dvt i\geq1\}.$ Let $ G_{n}$ be the sequence of empirical
distributions corresponding to $\{(x_{i},y_{i})\dvt i\geq1\}$ with common
distribution $G_{0}$. Then%
\begin{equation}
n^{1/2}\bigl(T_{\mathrm{MM},\beta}(G_{n})-\beta_{0}\bigr)\rightarrow_{d}N(0,V)\vspace*{-2pt}
\end{equation}
with\vspace*{-2pt}
\begin{equation}
V=\sigma_{0}^{2}\frac{\mathrm{E}_{F_{0}}\psi_{1}^{2} (  \fracb
{u-\alpha_{01}}{\sigma_{0}} )  }{\mathrm{E}_{F_{0}}^{2}\psi_{1}^{\prime
} (  \fracb{u-\alpha_{01}}{\sigma_{0}} )  }A_{0}^{-1}. \label{varerg}\vspace*{-2pt}
\end{equation}
A similar result can be obtained for
$T_{\mathrm{S},\beta}.$\vadjust{\goodbreak}

\item[(ix)] Assume that the conditions in \textup{(iv)} hold, that $\{(u_{i}%
,x_{i})\dvt i\geq1\}$ is  $\phi$-mixing (see, e.g., Billingsley
\cite{Billingsley1} for the definition of $\phi$-mixing)  with $\sum
_{i=1}^{\infty}\phi_{n}^{1/2}<\infty,$ that $u_{i}$ have a symmetric
distribution $F_{0}$ and that $\{x_{i}, i\geq1\}$  and $\{u_{i}\dvt i\geq1\}$
are independent. Let $ G_{n}$ be the sequence of empirical distributions
corresponding to $\{(x_{i},y_{i})\dvt i\geq1\}$ with common distribution $G_{0}$.
Then%
\begin{equation}
n^{1/2}\bigl(T_{\mathrm{MM}}(G_{n})-T_{\mathrm{MM}}(G_{0})\bigr) \rightarrow_{d}N(0,V),
\end{equation}
where
\begin{eqnarray} \label{varrmix}%
V&=&\frac{\sigma_{0}^{2} }{\mathrm{E}_{F_{0}}^{2}\psi_{1}^{\prime
} (  (u-\alpha_{0})/\sigma_{0} )  }C_{0}^{-1} \Biggl(
{\displaystyle\sum\limits_{i=-\infty}^{\infty}}
c_{i}C_{i} \Biggr)  C_{0}^{-1},
\nonumber\\[-2pt]
c_{i}&=&\mathrm{E} \biggl[  \psi_{1} \biggl(  \frac{u_{1}-\alpha_{0}}{\sigma_{0}%
} \biggr)  \psi_{1} \biggl(  \frac{u_{1+i}-\alpha_{0}}{\sigma_{0}} \biggr)
 \biggr]  ,
\\[-2pt]
C_{i}&=&\mathrm{E}\underline{\dot{g}}(x_{1},\beta_{0})\underline{\dot{g}%
}(x_{1+i},\beta_{0})^{\prime}\nonumber%
\end{eqnarray}
and $T_{\mathrm{MM}}(G_{0})=(\beta_{0},a_{0}).$
\end{enumerate}
\end{theorem}

\begin{remark}
Note that (viii) implies that the asymptotic covariance matrix of\break
$ n^{1/2}(T_{\mathrm{MM},\beta}(G_{n})-\beta_{0})$ is the same as when the $u_{i}$ are
i.i.d. This result does not hold for the intercept estimate
$T_{\mathrm{MM},\alpha}(G_{n})$. Croux, Dhaene and Hoorelbeke~\cite{Croux2} derived a
similar result for linear regression through the origin with one covariable
with mean 0.
\end{remark}

\begin{remark}
The $\phi$-mixing condition in (ix) can be replaced by any other type of
mixing condition  that guarantees the validity of  the central limit theorem
(see, e.g., Section 1.5.1 of  Doukhan~\cite{Doukham}). A result
similar to part (ix) of Theorem~\ref{teoremaso} was stated by Croux \textit{et al}.~\cite{Croux2}.
\end{remark}

\subsection{Location case}

The location model corresponds to the case where there are no regressors:
$p=q=0$ and so $y_{i}=u_{i}$ and $\xi=\alpha$. If $F_{0}$ denotes the common
distribution of the $u_{i}$, then $T(F_{0})=(T_{\mathrm{S}}(F_{0}%
),T_{\mathrm{MM}}(F_{0}),S(F_{0}))$ is defined  as in the regression case
with $\underline{g}(x,\xi)$  replaced by $\alpha.$ Then, the resulting
$T_{\mathrm{MM}}=T_{\mathrm{MM},\alpha}$ and $T_{\mathrm{S}}=T_{\mathrm{S}%
,\alpha}$ are the location functionals while $S$ is a functional estimating
the error scale. In this case, $I_{T_{\mathrm{MM}},F_{0}}$ is given by%
\begin{eqnarray}\label{INFMM3}%
I_{T_{\mathrm{MM}},F_{0}}(x)  &  =&\frac{\sigma_{0}}{a_{01}}\psi_{1} \biggl(
\frac{y-\alpha_{01}}{\sigma_{0}} \biggr) \nonumber
\\[-8pt]
\\[-8pt]
&&{}  -\frac{e_{01}\sigma_{0}}{a_{01}d_{0}}\biggl (  \rho_{0}\biggl (  \frac
{y-\alpha_{00}}{\sigma_{0}} \biggr)  -\delta \biggr)  .
\nonumber
\end{eqnarray}

The following theorem summarizes the properties of $T$ that can be derived
from the theorems in the former sections.

\begin{theorem}
\label{TeoLtotal} Assume that $\rho_{0}$ and $\rho_{1}$ are bounded $\rho
$-functions that satisfy Condition~\ref{R1}, with $\rho_{1}\leq\rho_{0}$. We
assume that $ F_{0}$ has a strong unimodal density. Then
\begin{enumerate}[(vii)]
\item[(i)] $T_{\mathrm{S}}(F_{0})=\alpha_{00}$ is the unique minimizer in
(\ref{S}). If $F_{0}$ is symmetric with respect to $\mu_{0}$, we have
$\alpha_{00}=\mu_{0.}$

\item[(ii)] $T_{\mathrm{\mathrm{MM}}}(F_{0})=\alpha_{01}$ is the unique
minimizer in (\ref{M-loc}).  If $F_{0}$ is symmetric with respect to $\mu
_{0}$, we have $\alpha_{01}=\mu_{0}$.

\item[(iii)] The functional $T=(T_{\mathrm{S}},T_{\mathrm{MM}},S)$ is weakly
continuous at $F_{0}.$

\item[(iv)] Assume also that Condition~\ref{R2}  holds and that $a_{00}%
\neq0,$ $a_{01}\neq0,$ $d_{0}\neq0$. Then, $D_{0}=\mathrm{E}_{F_{0}}\dot{\Psi
} (  z,T(F_{0}) )  $ is invertible, $I_{T_{\mathrm{MM}},F_{0}}%
(y)$ is given by (\ref{INFMM3}). The influence function $I_{T_{\mathrm{S}%
},F_{0}}(y)$  has a similar expression replacing  $\alpha_{01},a_{01}$ and
$e_{01}$ by $\alpha_{00},a_{00}$ and $e_{00}$, respectively.

\item[(v)] Under the same assumptions as in \textup{(iv)}, let $\{F_{n}\}$ be a
sequence of random distributions converging weakly to $F_{0}$ and satisfying
Condition~\ref{CondiC} a.s. Then $T$ is a.s. weakly differentiable at
$\{F_{n}\}$.

\item[(vi)] Assume the same conditions as in \textup{(v)} and%
\begin{equation}
n^{1/2}\mathrm{E}_{F_{n}}I_{T,F_{0}}(y)\rightarrow_{d}H. \label{convHL}
\end{equation}
Then
\begin{equation}
n^{1/2}\bigl(T(F_{n})-T(F_{0})\bigr)=n_{F_{n}}^{1/2}\mathrm{E}I_{T,F_{0}}(y)+\mathrm{o}_{p}(1),
\label{asinMML}%
\end{equation}
and therefore
\begin{equation}
n^{1/2}\bigl(T(F_{n})-T(F_{0})\bigr)\rightarrow_{d}H. \label{iv1L}
\end{equation}

\item[(vii)] Assume the same conditions as in \textup{(iv)}. Let $\{F_{n}\}$ be the
sequence of empirical distributions corresponding to i.i.d. observations
$u_{i}$ with common distribution $F_{0}$. Then (\ref{convHL}) holds with $H =N(0,V)$ and $V$ given by (\ref{defV}). If $F_{0}$ is symmetric,  the
asymptotic variance of $T_{\mathrm{MM}}$ given by (\ref{varsym}) becomes%
\[
V=\sigma_{0}^{2}\frac{\mathrm{E}_{F_{0}}\psi_{1} (  u/\sigma_{0} )
^{2}}{ (  \mathrm{E}_{F_{0}}\psi_{1}^{\prime} (  u/\sigma_{0} )
 )  ^{2}}.\vspace*{-4pt}
\]
\end{enumerate}
\end{theorem}

\section{Proofs}\label{SecProofs}\vspace*{-2pt}

\subsection{\texorpdfstring{Proof of Theorem \protect\ref{Teore5}}{Proof of Theorem 1}}\vspace*{-2pt}

We shall need the following auxiliary result, which is due to Ibragimov
\cite{Ibragimov}.\vspace*{-2pt}

\begin{lemma}
\label{Teore3}If $f$ is  a strongly unimodal density and $\varphi$ is a
density such that $\log\varphi$ is concave on its support, the convolution%
\begin{equation}
h(t)=\int_{-\infty}^{\infty}\varphi(u-t)f(u)\,\mathrm{d}u \label{hdet}\vadjust{\goodbreak}
\end{equation}
is strongly unimodal.\vadjust{\goodbreak}
\end{lemma}

The following lemma is a small variation of one given by Mizera~\cite{Mizera}.

\begin{lemma}
\label{Teore4} Let $\rho$ satisfy Condition~\ref{R1} and let $F$ be a
distribution with a strongly unimodal density~$f$. Then \textup{(a)} there exists
$t_{0}$ such that
\begin{equation}
q(t)=\mathrm{E}_{F}\rho(u-t) \label{qQ0}%
\end{equation}
 has a unique minimum at $t_{0};$ \textup{(b)} if $F$ is symmetric around $\mu_{0},$
then $t_{0}=\mu_{0}.$
\end{lemma}

\begin{pf}
(a) Put $k=\int_{-m}^{m}\rho(x)\,\mathrm{d}x$ and $\varphi(u)=(1-\rho(u))/k,$ which
vanishes for $|u|>m.$ Then%
\[
q(t)=1-\mathrm{E}_{F}\bigl(1-\rho(u-t)\bigr)=1-k\mathrm{E}_{F}\varphi(u-t)=1-kh(t),
\]
where $ h(t)$ is given by (\ref{hdet}). Since by Lemma~\ref{Teore3} $h(t)$
is a strongly unimodal density, part (a) of the lemma follows.

(b) It is proved in Lemma 3.1 of Yohai~\cite{YohaiI}.
\end{pf}

\begin{pf*}{Proof of Theorem~\ref{Teore5}}
Without loss of generality we may assume $\sigma=1.$ To prove the theorem, we
will show that the unique minimum of $ R(\beta,\alpha)=\mathrm{E}_{G_{0}%
}\rho(y-g(x,\beta)-\alpha)$ is $\beta=\beta_{0},\alpha=t_{0}.$  We will first
prove that%
\[
R(\beta_{0},t_{0})<R(\beta_{0},\alpha) \qquad \mbox{for }\alpha\neq t_{0}.
\]
This is equivalent to
\[
\mathrm{E}_{F_{0}}\rho(u-t_{0})<\mathrm{E}_{F_{0}}\rho(u-\alpha) \qquad \mbox{for
}\alpha\neq t_{0},
\]
 which follows from Theorem~\ref{Teore4}.

Consider now $(\beta,\alpha)$ with $\beta\neq\beta_{0}.$ Let $A=\{x\dvt g(x,\beta
_{0})=g(x,\beta)+\alpha-t_{0}\}$  and $q$ as in (\ref{qQ0}), with $F$
replaced by $F_{0}$. Then%
\begin{eqnarray}\label{espesp}
R(\beta,\alpha)  &  =&\mathrm{E}_{G_{0}}\bigl\{\mathrm{E}_{G_{0}}\bigl[\rho
\bigl(y-g(x,\beta)-\alpha\bigr)|x\bigr]\bigr\}\nonumber
\\[-8pt]
\\[-8pt]
&  =&\mathrm{E}_{G_{0}}\bigl\{\mathrm{E}_{G_{0}}\bigl[\rho\bigl(u+g(x,\beta_{0})-g(x,\beta
)-\alpha\bigr)|x\bigr]\bigr\}. %
\nonumber
\end{eqnarray}
Since $ u$ and $x$ are independent, we get%
\begin{equation}
\mathrm{E}\bigl[\rho\bigl(u+g(x,\beta_{0})-g(x,\beta)-\alpha\bigr)|x\bigr]=q\bigl(g(x,\beta
)-g(x,\beta_{0})+\alpha\bigr). \label{Econdi}%
\end{equation}
Then according to Theorem~\ref{Teore4}, the left-hand side of (\ref{Econdi})
is equal to $q(t_{0})$ if $x\in A$ and grater than $q(t_{0})$ otherwise.
Condition~\ref{CondId} implies that $\mathrm{P}(A^{c})>0$ and from
(\ref{espesp}) we get that $R(\beta,\alpha)>q(t_{0}).$ Finally, the theorem
follows from the fact that $R(\beta_{0},t_{0})=q(t_{0}).$
\end{pf*}

\subsection{\texorpdfstring{Proof of Theorems \protect\ref{Teore01} and \protect\ref{Teore02}}
{Proof of Theorems 2 and 3}}

Before proving Theorems~\ref{Teore01} and~\ref{Teore02}, we need some auxiliary
results.\vadjust{\goodbreak}

\begin{lemma}
\label{Lema1} Consider distributions $\{G_{n}\}$ and $G_{0}$ on
$R^{p}\times R$. Let $\{\xi_{n}\}$ and $\{\sigma_{n}\}$ be sequences in
$B\times R$ and $R_{+}$, respectively, such that $\xi_{n}\rightarrow\xi\in
B\times R$ and $\sigma_{n}\rightarrow\sigma>0$. Assume that $\underline{g}( x,
\xi) $ is continuous in $\xi$. If $G_{n}\rightarrow_{w}G_{0}$, then
\[
\lim_{n\rightarrow\infty}\mathrm{E}_{G_{n}}\rho \biggl(  \frac
{y-\underline{g}(x,\xi_{n})}{\sigma_{n}} \biggr)  =\mathrm{E}_{G_{0}}%
\rho \biggl(  \frac{y-\underline{g}(x,\xi)}{\sigma} \biggr)  .
\]
\end{lemma}

\begin{pf}
Since $G_{n}\rightarrow_{w}G_{0}$ and $\rho$  is  continuous
and bounded, we have
\[
\mathrm{E}_{G_{n}}\rho \biggl(  \frac{y-\underline{g}(x,\xi)}{\sigma} \biggr)
\rightarrow\mathrm{E}_{G_{0}}\rho \biggl(  \frac{y-\underline{g}(x,\xi)}{\sigma
} \biggr)  ,
\]
and therefore it suffices to show that
\[
\mathrm{E}_{G_{n}}\rho \biggl(  \frac{y-\underline{g}(x,\xi_{n})}{\sigma_{n}%
} \biggr)  -\mathrm{E}_{G_{n}}\rho \biggl(  \frac{y-\underline{g}(x,\xi)}{\sigma
} \biggr)  \rightarrow0.
\]
Since $\{G_{n}\}_{n\geq1}$ is tight, it suffices to show that
if ${\mathcal{P}}$ is a tight set of distributions of $(x,y)$, then
\[
\sup_{F\in{\mathcal{P}}} \biggl| \mathrm{E}_{F}\rho\biggl (  \frac
{y-\underline{g}(x,\xi_{n})}{\sigma_{n}} \biggr)  -\mathrm{E}_{F}\rho \biggl(
\frac{y-\underline{g}(x,\xi)}{\sigma} \biggr)   \biggr| \rightarrow0.
\]

To prove this, put $z=(x,y)$. Then for all $K>0$%
\begin{eqnarray}\label{conn0}
&&   \biggl| \mathrm{E}_{F}\rho \biggl(  \frac{y-\underline{g}(x,\xi_{n}%
)}{\sigma_{n}} \biggr)  -\mathrm{E}_{F}\rho \biggl(  \frac{y-\underline{g}%
(x,\xi)}{\sigma} \biggr)   \biggr| \nonumber
\\[-8pt]
\\[-8pt]
&& \quad   \leq2\mathrm{E}_{F}\mathbf{1}_{\{ \|z\|>K\}}+\mathrm{E}_{F}
\biggl|
\rho \biggl(  \frac{y-\underline{g}(x,\xi_{n})}{\sigma_{n}} \biggr)  -\rho \biggl(
\frac{y-\underline{g}(x,\xi)}{\sigma} \biggr)   \biggr| \mathbf{1}%
_{\{ \|z\|\leq K\}}.
\nonumber
\end{eqnarray}

If $\|z\|\leq K$, we have
\begin{eqnarray}\label{conn}
& &  \biggl| \frac{y-\underline{g}(x,\xi_{n})}{\sigma_{n}}-\frac
{y-\underline{g}(x,\xi)}{\sigma} \biggr|\nonumber
\\[-8pt]
\\[-8pt]
&& \quad   \leq\frac{1}{\sigma\sigma_{n}}[|\sigma_{n}-\sigma||y|+|\sigma_{n}%
-\sigma||\underline{g}(x,\xi)|+\sigma|\underline{g}(x,\xi_{n})-\underline
{g}(x,\xi)|]. %
\nonumber
\end{eqnarray}
Now, given $\varepsilon>0$,  we can find $K$ such that
\[
2\sup_{F\in{\mathcal{P}}}\mathrm{P}_{F} (  \|z\|>K )  \leq
\varepsilon/2
\]
and $\alpha$ such that
\[
|\rho(u)-\rho(v)|\leq\varepsilon/2 \qquad \mbox{if }|u-v|\leq\alpha.
\]
Then, we can choose $n_{0}$ such that the right-hand side of (\ref{conn}) is
smaller than $\alpha$ if $n\geq n_{0}$ and $\|z\|\leq K$, and so from
(\ref{conn0})  we obtain for all $n\geq n_{0}$
\[
 \biggl| \mathrm{E}_{F}\rho \biggl(  \frac{y-\underline{g}(x,\xi_{n})}%
{\sigma_{n}} \biggr)  -\mathrm{E}_{F}\rho \biggl(  \frac{y-\underline{g}(x,\xi
)}{\sigma} \biggr)   \biggr|  \leq \varepsilon \qquad  \forall F\in{\mathcal{P}%
}.
\]
\upqed
\end{pf}

\begin{lemma}
\label{Lema2}Assume that $B$ is closed  and let $G_{0}$ be any distribution
for $(x,y)$ such that (\ref{M-loc}) has a unique solution $T_{\mathrm{M}%
}(G_{0})$. Let $\{G_{n}\}$ be a sequence such that $G_{n}\rightarrow_{w}G_{0}$
and $\{T_{\mathrm{M}}(G_{n})\}$ is bounded. If $\widetilde{S}(G_{n}%
)\rightarrow\widetilde{S}(G_{0})>0$, then $T_{\mathrm{M}}(G_{n})\rightarrow
T_{\mathrm{M}}(G_{0})$.
\end{lemma}

\begin{pf} Put for brevity
\begin{equation}
\xi_{n}=T_{\mathrm{M}}(G_{n}), \qquad  \xi_{0}=T_{\mathrm{M}}(G_{0}), \qquad  \sigma
_{n}=\widetilde{S}(G_{n}), \qquad  \sigma_{0}=\widetilde{S}(G_{0}). \label{breve}%
\end{equation}
Since $\{\xi_{n}\}$ remains in a compact set, it suffices to prove that
$\xi_{0}$ is the only accumulation point of $\{\xi_{n}\}$, that is, if a
subsequence tends to some $\widehat{\xi}$, then $\widehat{\xi}=\xi_{0}$.
Without loss of generality, assume that $\xi_{n}\rightarrow\widehat{\xi}$. The
definition of $\xi_{n}$ implies
\begin{equation}
\mathrm{E}_{G_{n}}\rho \biggl(  \frac{y-\underline{g}(x,\xi_{n})}{\sigma_{n}%
} \biggr)  \leq\mathrm{E}_{G_{n}}\rho \biggl(  \frac{y-\underline{g}(x,\xi_{0}%
)}{\sigma_{n}} \biggr)  . \label{desigualdad}%
\end{equation}
Using Lemma~\ref{Lema1}, we get%
\[
M_{G_{0}}(\widehat{\xi})=\mathrm{E}_{G_{0}}\rho \biggl(  \frac{y-\underline
{g}(x,\widehat{\xi})}{\sigma_{0}} \biggr)  \leq\mathrm{E}_{G_{0}}\rho\biggl (
\frac{y-\underline{g}(x,\xi_{0})}{\sigma_{0}} \biggr)  =M_{G_{0}}(\xi_{0}).
\]
Since $\xi_{0}$ is the only minimizer of $M_{G_{0}}$, we conclude that
$\widehat{\xi}=\xi_{0}$.
\end{pf}

\begin{lemma}
\label{Lema3}
Let $\{\xi_{n}\}$ and $\{\sigma_{n}\}$ be sequences in $R^{p+1}$
and $R_{+},$ respectively. Assume that when $n\rightarrow\infty,$
$G_{n}\rightarrow_{w}G_{0}$, $\|\xi_{n}\|\rightarrow\infty$ and $\{\sigma
_{n}\}$ is bounded. Then
\begin{equation}
\lim\inf_{n\rightarrow\infty}\mathrm{E}_{G_{n}}\rho \biggl(  \frac
{y-\xi_{n}^{\prime}(x^{\prime},1)^{\prime}}{\sigma_{n}} \biggr)  \geq1-c_{0},
\label{cc}%
\end{equation}
where $c_{0}=c (  G_{0} )  $ is defined in (\ref{c}).
\end{lemma}

\begin{pf} Assume without loss of generality that there exist $\gamma\in
R^{p}$ and $\sigma>0$  such that for some subsequence $\gamma_{n}=\xi
_{n}/\|\xi_{n}\|\rightarrow\gamma$, and $ \sigma_{n}\leq\sigma$. Put
$\lambda_{n}=\|\xi_{n}\|.$

For $\varepsilon>0$ let $d_{\varepsilon}$ be such that $\rho(u)\geq
1-\varepsilon$ for $|u|\geq d_{\varepsilon}$. Therefore,
\[
\mathrm{E}_{G_{n}}\rho \biggl(  \frac{y-\xi_{n}^{\prime}(x^{\prime},1)^{\prime}%
}{\sigma_{n}} \biggr)  \geq\mathrm{E}_{G_{n}}\rho \biggl(  \frac{y-\xi
_{n}^{\prime}(x^{\prime},1)^{\prime}}{\sigma} \biggr)  \geq(1-\varepsilon
)\mathrm{P}_{G_{n}} \biggl(  \frac{|y-\lambda_{n}\boldsymbol{\gamma}_{n}^{\prime
}(x^{\prime},1)^{\prime}|}{\sigma}\geq d_{\varepsilon} \biggr)  .
\]

Then, to prove the lemma, it suffices to show that
\[
\lim\inf_{n\rightarrow\infty}\mathrm{P}_{G_{n}} \biggl(   \biggl|
\frac{y}{\lambda_{n}}-\boldsymbol{\gamma}_{n}^{\prime}(x^{\prime},1)^{\prime
} \biggr| \geq\frac{d_{\varepsilon}\sigma}{\lambda_{n}} \biggr)  \geq
1-c_{0}.
\]

Let $(x_{n},y_{n})\sim G_{n}$ and $(x_{0},y_{0})\sim G_{0}$. Since
$\lambda_{n}\rightarrow\infty$, we have $y_{n}/\lambda_{n}\rightarrow_{p}0$.
Then the convergence of $\gamma_{n}$ to $\gamma$ guarantees that
\[
\frac{y_{n}}{\lambda_{n}}-\boldsymbol{\gamma}_{n}^{\prime}(x_{n}^{\prime
},1)^{\prime}\rightarrow_{d}\gamma^{\prime}(x_{0}^{\prime},1)^{\prime}.
\]

For any $\alpha>0$ which is a point of continuity of the distribution of
$|\gamma^{\prime}(x_{0},1)|$, $\lambda_{n}\rightarrow\infty$ implies%
\begin{eqnarray*}
\lim\inf_{n\rightarrow\infty}\mathrm{P}_{G_{n}} \biggl(   \biggl|
\frac{y}{\lambda_{n}}-\boldsymbol{\gamma}_{n}^{\prime}(x^{\prime},1)^{\prime
} \biggr| >\frac{d_{\varepsilon}\sigma}{\lambda_{n}} \biggr)   &\geq&\lim
\inf_{n\rightarrow\infty}\mathrm{P}_{G_{n}} \biggl(   \biggl|\frac
{y}{\lambda_{n}}-\boldsymbol{\gamma}_{n}^{\prime}(x^{\prime},1)^{\prime
} \biggr|>\alpha \biggr)\\
  &=&\mathrm{P}_{G_{0}}\bigl (  |\gamma^{\prime}(x^{\prime
},1)^{\prime}|>\alpha \bigr)  .
\end{eqnarray*}
Letting $\alpha\rightarrow0$ and recalling (\ref{c}), we get%
\[
\lim\inf_{n\rightarrow\infty}\mathrm{P}_{G_{n}} \biggl(   \biggl|\frac
{y}{\lambda_{n}}-\boldsymbol{\gamma}_{n}^{\prime}(x^{\prime},1)^{\prime
} \biggr|>\frac{d_{\varepsilon}\sigma}{\lambda_{n}} \biggr)  \geq1-c_{0}.
\]
\upqed
\end{pf}

The proof of the following lemma is similar to that of Lemma~\ref{Lema3}.

\begin{lemma}
\label{Lema4}
Let $\{\xi_{n}\}$ be a sequence in $B\times R$, with $B$ compact.
Assume that when $n\rightarrow\infty,$ $G_{n}\rightarrow_{w}G_{0}$, $\|\xi
_{n}\|\rightarrow\infty$ and $\{\sigma_{n}\}$ is bounded. Then
\begin{equation}
\lim\inf_{n\rightarrow\infty}\mathrm{E}_{G_{n}}\rho\biggl (
\frac{y-\underline{g}(x,\xi_{n})}{\sigma_{n}} \biggr)  =1. \label{Lemma4eq}%
\end{equation}
\end{lemma}

Finally, the following result we be used.

\begin{lemma}
\label{s1} Let $S(G)$ be defined by (\ref{defSG}) and suppose that
 $S(G_{0})>0.$ Then,  $G_{n}\rightarrow_{w}G_{0}$ implies that there exists
$n_{0}$  such that $S(G_{n})>0$  for $n\geq n_{0}$.
\end{lemma}

\begin{pf} Suppose that the lemma is not true. Then there exists a
subsequence $\{G_{n_{k}}\}_{k\geq1}$ such that $S(G_{n_{k}})=0$ for all $k$.
This means that giving $\varepsilon>0$, there exists $(\beta_{n_{k}}%
,\alpha_{n_{k}})$ such that
\[
\mathrm{E}_{G_{n_{k}}}\rho_{0} \biggl(  \frac{y-g(\mathbf{x},\beta_{n_{k}%
})-\alpha_{n_{k}}}{\varepsilon} \biggr)  <\delta  \qquad \mbox{for any $s>0$}.
\]
The same arguments that we use to prove Lemma~\ref{Lema3} let us show that
$\{(\beta_{n_{k}},\alpha_{n_{k}})\}$ is bounded and  therefore (passing on to
a subsequence if necessary) we can assume that $(\beta_{n_{k}},\alpha_{n_{k}%
})\rightarrow(\tilde{\beta},\tilde{\alpha})$. Then, from Lemma~\ref{Lema1} we
get that
\[
\mathrm{E}_{G_{0}}\rho_{0} \biggl(  \frac{y-g(\mathbf{x},\tilde{\beta}%
)-\tilde{\alpha}}{\varepsilon} \biggr)  \leq\delta  \qquad \mbox{for any $s>0$}.
\]

Then, $S(G_{0})\leq S^{\ast}(G_{0},\tilde{\beta},\tilde{\alpha})\leq
\varepsilon$. Since this holds for any $\varepsilon>0$,  we get that
$S(G_{0})=0$. This contradicts the assumption that $S(G_{0})>0.$
\end{pf}

\subsubsection{\texorpdfstring{Proof of Theorem \protect\ref{Teore01}}{Proof of Theorem 2}}

Let $G_{n}\rightarrow_{w}G_{0}$. Since $\widetilde{S}$ is weakly continuous at
$G_{0}$, it follows that $\widetilde{S}(G_{n})\rightarrow\widetilde{S}%
(G_{0})>0$, by hypothesis.

Case (a):  We prove first that $\{T_{\mathrm{M}}(G_{n})\}$ is bounded.
 Suppose that it is not true; then without loss of generality we may assume
that $\|T_{\mathrm{M}}(G_{n})\|\rightarrow\infty.$ Then Lemma~\ref{Lema4} implies%
\[
1 = \liminf_{n\rightarrow\infty}M_{G_{n}} (  T_{\mathrm{M}}%
(G_{n}) )   \leq \liminf_{n\rightarrow\infty}M_{G_{n}} (
T_{\mathrm{M}}(G_{0}) )  =M_{G_{0}} (  T_{\mathrm{M}}(G_{0}) )
,
\]
and this implies that $M_{G_{0}}(\xi)=1$ for all $\xi$. This contradicts the
assumption  that $T_{\mathrm{M}}(G_{0})$ is univocally defined. Then, $\{T_{\mathrm{M}%
}(G_{n})\}$ is bounded and from Lemma~\ref{Lema2}, we get that $T_{\mathrm{M}%
}(G_{n})\to T_{\mathrm{M}}(G_{0})$.

Case (b): Recall the notation in (\ref{breve}). Convergence of $\{\sigma
_{n}\}$ guarantees that it is a bounded sequence. Suppose that $\{\xi_{n}\}$
is unbounded. Then, passing on to a subsequence if necessary, we may assume
that $\|\xi_{n}\|\rightarrow\infty.$ In this case by Lemma~\ref{Lema3} we have%
\begin{equation}
\lim\inf_{n\rightarrow\infty}M_{G_{n}}(\xi_{n})=\lim\inf
_{n\rightarrow\infty}\mathrm{E}_{G_{n}}\rho \biggl(  \frac{y-\xi
_{n}^{^{\prime}}(x^{\prime},1)^{\prime}}{\sigma_{n}} \biggr)  \geq1-c_{0}.
\label{cont10}%
\end{equation}

We also have%
\begin{equation}
\lim_{n\rightarrow\infty}M_{G_{n}}(\xi_{0})=\lim_{n\rightarrow
\infty}\mathrm{E}_{G_{n}}\rho \biggl(  \frac{y-\xi_{0}^{^{\prime}}(x^{\prime
},1)^{\prime}}{\sigma_{n}} \biggr)  =M_{G_{0}}(\xi_{0})<1-c_{0}. \label{cont20}%
\end{equation}
Inequalities (\ref{cont10}) and (\ref{cont20}) imply that  there exists
$n_{0}$ such that for $n\geq n_{0}$
\[
M_{G_{n}}(\xi_{n})>M_{G_{n}}(\xi_{0}),
\]
contradicting the definition of $T_{\mathrm{M}}(G_{n}).$ Therefore, $\{\xi
_{n}\}$ is bounded,  and  then the weak continuity of $T_{\mathrm{M}}$
follows from Lemma~\ref{Lema2}.

\subsubsection{\texorpdfstring{Proof of Theorem \protect\ref{Teore02}}{Proof of Theorem 3}}

Let $G_{n}\rightarrow_{w}G_{0}$, $\xi_{n}=T_{\mathrm{S}}(G_{n})$, $\xi
_{0}=T_{\mathrm{S}}(G_{0}),$ $\sigma_{n}=S(G_{n})$ and $\sigma_{0}=S(G_{0}).$
We prove first that $\{\sigma_{n}\}$ is bounded. Take  any $\sigma
_{1}>\sigma_{0};$ then by Lemma~\ref{Lema1}%
\[
\mathrm{E}_{G_{n}}\rho_{0} \biggl(  \frac{y-\underline{g}(x,\xi_{0})}{\sigma
_{1}} \biggr)  \rightarrow\mathrm{E}_{G_{0}}\rho_{0} \biggl(  \frac
{y-\underline{g}(x,\xi_{0})}{\sigma_{1}} \biggr)  <\delta,
\]
and therefore there exists $n_{0}$ such that
\begin{equation}
S^{\ast}(\xi_{0},G_{n})<\sigma_{1} \qquad \mbox{for }n\geq n_{0}, \label{nosec}%
\end{equation}
which implies that $S^{\ast}(G_{n},\xi_{0})$ is bounded and therefore
$\sigma_{n}\leq S^{\ast}(\xi_{0},G_{n})$ is also bounded. On the other hand,
by Lemma~\ref{s1}, we get that $\sigma_{n}>0$ for $n$ large enough.

We  now prove that  $\{\xi_{n}\}$ is bounded. In case (a), if $\{\xi_{n}\}$
is unbounded, Lemma~\ref{Lema4} implies%
\begin{equation}
\lim\inf_{n\rightarrow\infty}\mathrm{E}_{G_{n}}\rho_{0} \biggl(
\frac{y-\underline{g}(x,\xi_{n})}{\sigma_{n}} \biggr)  \geq1,
\end{equation}
and this contradicts the fact that for all $n$%
\[
\mathrm{E}_{G_{n}}\rho_{0} \biggl(  \frac{y-\underline{g}(x,\xi_{n})}{\sigma
_{n}} \biggr)  =\delta<1.
\]

Consider now case (b) and assume that $\{\xi_{n}\}$ is unbounded. Then,
passing on to a subsequence  if necessary,  we may assume that $\|\xi
_{n}\|\rightarrow\infty.$ Then by Lemma~\ref{Lema3}
\[
\lim\inf_{n\rightarrow\infty}\mathrm{E}_{G_{n}}\rho_{0}\biggl (
\frac{y-\xi_{n}^{\prime}(x^{\prime},1)^{\prime}}{\sigma_{n}} \biggr)
\geq1-c_{0},
\]
and this contradicts the fact that for all $n$%
\[
\mathrm{E}_{G_{n}}\rho_{0} \biggl(  \frac{y-\xi_{n}^{\prime}(x^{\prime
},1)^{\prime}}{\sigma_{n}} \biggr)  =\delta<1-c_{0}.
\]
Then in case (b) $\{\xi_{n}\}$ is also bounded.

We now show that $\sigma_{n}\rightarrow\sigma_{0}.$ Suppose that this is not
true.  By  passing on to a subsequence if necessary,  we may assume that
$\sigma_{n}\rightarrow\sigma^{\ast}\neq\sigma_{0}$ and $\xi_{n}\rightarrow
\xi^{\ast}$ for some $\xi^{\ast}$ and $\sigma^{\ast}.$ Since (\ref{nosec})
holds for any $\sigma^{\prime}>\sigma_{0}$ we have $\sigma^{\ast}\leq
\sigma_{0}$ and therefore $\sigma^{\ast}<\sigma_{0}.$ Then Lemma~\ref{Lema1}
implies%
\[
\delta=\lim_{n\rightarrow\infty}\mathrm{E}_{G_{n}}\rho_{0} \biggl(
\frac{y-\underline{g}(\xi_{n},x)}{\sigma_{n}} \biggr)  =\mathrm{E}_{G_{0}}%
\rho_{0} \biggl(  \frac{y-\underline{g}(\xi^{\ast},x)}{\sigma^{\ast}} \biggr)  ,
\]
and therefore $S (  G_{0} )  \leq S^{\ast}(G_{0},$$\xi^{\ast}%
)=\sigma^{\ast}<\sigma_{0}$. This contradicts the fact that $S(G_{0}%
)=\sigma_{0}$ and shows that $S$ is weakly continuous.

Finally, the weak continuity of $T_{\mathrm{S}}$ follows from  (\ref{s=m})
and Theorem~\ref{Teore01}.

\subsection{\texorpdfstring{Proofs of Theorems \protect\ref{TeoDiffe} and \protect\ref{Teorevarios}}
{Proofs of Theorems 4 and 5}}

\subsubsection{\texorpdfstring{Proof of Theorem \protect\ref{TeoDiffe}}{Proof of Theorem 4}}

Since%
\[
\mathrm{E}_{G_{n}}\Psi(z,T(G_{n}))=0,
\]
the Mean Value theorem together with Condition~\ref{CondiC} and the
consistency of $T (  G_{n} )  $ yield%
\[
\mathrm{E}_{G_{n}}\Psi(z,T(G_{0}))+D (  G_{n},\theta_{n}^{\ast
} )  \bigl(T(G_{n})-T(G_{0})\bigr)=0,
\]
where $\theta_{n}^{\ast}\rightarrow\theta_{0}.$ Then,  (\ref{cca}) implies
that $D (  G_{n},\theta_{n}^{\ast} )  \rightarrow D_{0}$ and,  since
for large $n$, $D (  G_{n},\theta_{n}^{\ast} )  $ is nonsingular, we
may write%
\begin{eqnarray*}
T(G_{n})-T(G_{0})  &  =&-D (  G_{n},\theta_{n}^{\ast} )
^{-1}\mathrm{E}_{G_{n}}\Psi(z,T(G_{0}))\\[-2pt]
&  =&\mathrm{E}_{G_{n}}I_{T,G_{0}}(z)+ \bigl(  D_{0}^{-1}-D (  G_{n}%
,\theta_{n}^{\ast} )  ^{-1} \bigr)  \mathrm{E}_{G_{n}}I_{T,G_{0}}(z).
\end{eqnarray*}
Condition~\ref{CondiA}  implies that the second term of the right--hand
side divided by $\|\mathrm{E}_{G_{n}}I_{T,G_{0}}(z)\|$ tends to zero, and
 this proves the theorem.

\subsubsection{\texorpdfstring{Proof of Theorem \protect\ref{Teorevarios}}{Proof of Theorem 5}}

Under the assumptions of this theorem, we can prove that Condition
\ref{CondiA} holds a.s. using  the same arguments as in Lemma 4.2 of Yohai
\cite{YohaiI}.  The only change is to replace the Law of Large Numbers for
i.i.d. random variables by the assumption  that $\mathrm{E}_{G_{n}%
}d(z)\rightarrow\mathrm{E}_{G_{0}}d(z)$ a.s. for all $d$ such that
$\mathrm{E}_{G_{0}}|d(z)|<\infty$ in the case (a) and for the fact that
$\mathrm{E}_{G_{n}}d(z)\rightarrow\mathrm{E}_{G_{0}}d(z)$ for all function $d$
bounded and continuous in case (b). Then, Theorem~\ref{TeoDiffe} implies that
$T$ is weakly differentiable at $\{G_{n}\}$.

\subsection{Derivations of influence functions}\label{secDerIF}

\subsubsection{\texorpdfstring{Derivation of (\protect\ref{INFMM1})--(\protect\ref{INFMM2})}
{Derivation of (6.11)--(6.12)}}

Put for brevity%
\[
t_{\mathrm{MM}}=\frac{y-\underline{g}(x,\xi_{\mathrm{MM}})}{\sigma
}, \qquad  t_{\mathrm{S}}=\frac{y-\underline{g}(x,\xi_{\mathrm{S}})}{\sigma}.
\]
Then
\[
\dot{\Psi} (  z,\theta )  = \left[
\begin{array}{c@{ \quad }c@{ \quad }c}%
\dot{\Psi}_{11} (  z,\theta )  & 0 & \dot{\Psi}_{13} (
z,\theta ) \\[3pt]
0 & \dot{\Psi}_{22} (  z,\theta )  & \dot{\Psi}_{23} (
z,\theta ) \\[3pt]
\dot{\Psi}_{31} (  z,\theta )  & 0 & \dot{\Psi}_{33} (
z,\theta )
\end{array}
 \right]  ,
\]
where
\begin{eqnarray}\label{matder1}
\dot{\Psi}_{11} (  z,\theta )   &  =&-\frac{1}{\sigma}\psi_{0}%
^{\prime} (  t_{\mathrm{S}} )  \underline{\dot{g}}(x,\xi_{\mathrm{S}%
})\underline{\dot{g}}(x,\xi_{\mathrm{S}})^{\prime}+\psi_{0} (
t_{\mathrm{S}} )  \underline{\ddot{g}} (  x,\xi_{\mathrm{S}} ),
\nonumber\\
\dot{\Psi}_{13} (  z,\theta )   &  =&-\frac{1}{\sigma}\psi_{0}%
^{\prime} (  t_{\mathrm{S}} )  t_{\mathrm{S}}\underline{\dot{g}%
}(x,\xi_{\mathrm{S}}),\nonumber\\
\dot{\Psi}_{22} (  z,\theta )   &  =&-\frac{1}{\sigma}\psi_{1}%
^{\prime} (  t_{\mathrm{MM}} )  \underline{\dot{g}}(x,{\xi
}_{\mathrm{MM}})\underline{\dot{g}}(x,\xi_{\mathrm{MM}})^{\prime}+\psi
_{1} (  t_{\mathrm{MM}} )  \underline{\ddot{g}} (  x,\xi
_{\mathrm{MM}} ), \nonumber
\\[-8pt]
\\[-8pt]
\dot{\Psi}_{23} (  z,\theta )   &  =&-\frac{1}{\sigma}\psi_{1}%
^{\prime} (  t_{\mathrm{MM}} )  t_{\mathrm{MM}}\underline{\dot{g}%
}(x,\xi_{\mathrm{MM}}),\nonumber\\
\dot{\Psi}_{31} (  z,\theta )   &  =&-\frac{1}{\sigma}\psi_{0} (
t_{\mathrm{S}} )  \underline{\dot{g}}(x,\xi_{\mathrm{S}}),\nonumber\\
\dot{\Psi}_{33} (  z,\theta )   &  =&-\frac{1}{\sigma}\psi_{0} (
t_{\mathrm{S}} )  t_{\mathrm{S}}.\nonumber\vadjust{\goodbreak}
\end{eqnarray}

From (\ref{matder1}) it is easy to show that
\[
D_{0}=\mathrm{E}_{G_{0}}\dot{\Psi}(z,\theta_{0})=-\frac{1}{\sigma_{0}} \left[
\begin{array}
[c]{c@{ \quad }c@{ \quad }c}%
a_{00}C_{0} & 0 & e_{00}b_{0}^{\ast}\\[3pt]
0 & a_{01}C_{0} & e_{01}b_{0}^{\ast}\\[3pt]
0 & 0 & d_{0}%
\end{array}
 \right]  .
\]
Therefore, $|D_{0}|=a_{00}a_{01}d_{0}|C_{0}|^{2}.$ It follows from
(\ref{defC0}) $|C_{0}|\neq0$ if and on only if $|A_{0}|\neq0,$ and that%
\[
C_{0}^{-1}= \left[
\begin{array}
[c]{c@{ \quad }c}%
A_{0}^{-1} & -A_{0}^{-1}b_{0}\\[2pt]
-(A_{0}^{-1}b_{0})^{\prime} & 1+b_{0}^{\prime}A_{0}^{-1}b_{0}%
\end{array}
 \right]  ,
\]
Direct calculation shows that%
\[
D_{0}^{-1}=-\sigma_{0}\left [
\begin{array}
[c]{c@{ \quad }c@{ \quad }c}%
a_{00}^{-1}C_{0}^{-1} & 0 & -e_{00}a_{00}^{-1}d_{0}^{-1}C_{0}^{-1}b_{0}^{\ast
}\\[3pt]
0 & a_{01}^{-1}C_{0}^{-1} & -e_{01}a_{01}^{-1}d_{0}^{-1}C_{0}^{-1}b_{0}^{\ast
}\\[3pt]
0 & 0 & d_{0}^{-1}%
\end{array}
 \right]  ,
\]
and the desired results follow from (\ref{inffund}).

\subsubsection{\texorpdfstring{Derivation of (\protect\ref{INFMM3})}{Derivation of (6.22)}}

In this case from (\ref{matder1}), it is easy to show that
\[
D_{0}=-\frac{1}{\sigma_{0}} \left[
\begin{array}
[c]{c@{ \quad }c@{ \quad }c}%
a_{00} & 0 & e_{00}\\
0 & a_{01} & e_{01}\\
0 & 0 & d_{0}%
\end{array}
 \right]  ,
\]
which implies%
\[
D_{0}^{-1}=-\sigma_{0} \left[
\begin{array}
[c]{c@{ \quad }c@{ \quad }c}%
a_{00}^{-1} & 0 & -e_{00}a_{00}^{-1}d_{0}^{-1} \\[3pt]
0 & a_{01}^{-1} & -e_{01}a_{01}^{-1}d_{0}^{-1} \\[3pt]
0 & 0 & d_{0}^{-1}%
\end{array}
 \right]  .
\]
The rest of the derivation is straightforward.

\subsection{\texorpdfstring{Proof of Theorems \protect\ref{teoremaso} and \protect\ref{TeoLtotal}}
{Proof of Theorems 6 and 7}}

\subsubsection{\texorpdfstring{Proof of Theorem \protect\ref{teoremaso}}{Proof of Theorem 6}}

Parts (i) and (ii) follow from Theorem~\ref{Teore5} and  Remark~\ref{Rt5}. To
prove (iii), we need to check conditions of Theorem~\ref{Teore01} and Theorem
\ref{Teore02}.  We start showing that $S(G_{0})>0.$ Let
\[
h_{\beta,\alpha}(s)=\mathrm{E}\rho_{0}\biggl (  \frac{y_{i}-g(x_{i}%
,\beta)-\alpha}{s} \biggr)  .
\]
Then, we have
\begin{equation}
\lim_{s\rightarrow\infty}h_{\beta,\alpha}(s)=\rho_{0}(0)=0 \label{lim1}
\end{equation}
and
\begin{equation}
\lim_{s\rightarrow0}h_{\beta,\alpha}(s)=1-\mathrm{P}\bigl(y_{i}=g(x_{i}%
,\beta)+\alpha\bigr). \label{lim2}%
\end{equation}
Since $ u_{i}$  has a continuous distribution and is independent of $x_{i}$,
we also have%
\begin{eqnarray}
\mathrm{P}\bigl(y_{i}=g(x_{i},\beta)+\alpha\bigr)&=&\mathrm{P}\bigl(g(x_{i},\beta_{0}%
)+u_{i}=g(x_{i},\beta)+\alpha\bigr)\nonumber
\\[-8pt]
\\[-8pt]
&=&\mathrm{E} \bigl[  \mathrm{P}\bigl(u_{i}%
=g(x_{i},\beta)-g(x_{i},\beta_{0})+\alpha\bigr) \bigr]  =0. \label{lim3}%
\nonumber
\end{eqnarray}
Equations (\ref{lim1}), (\ref{lim2}) and (\ref{lim3}) imply that $S^{\ast
}(G_{0},\beta,\alpha)>0$ for all $(\beta,\alpha)$, and so $S(G_{0})=S^{\ast
}(G_{0},\beta_{0},\alpha_{01})>0$.

Note that
\begin{eqnarray*}
M_{G_{0}}(T_{\mathrm{MM}}(G_{0}))  &  =&\mathrm{E} \biggl(  \rho_{1} \biggl(
\frac{y-T_{\mathrm{MM}}(G_{0})}{S(G_{0})} \biggr)   \biggr) \\
&  \leq&\mathrm{E} \biggl(  \rho_{1} \biggl(  \frac{y-T_{\mathrm{S}}(G_{0}%
)}{S(G_{0})} \biggr)   \biggr)   \leq\mathrm{E} \biggl(  \rho_{0} \biggl(  \frac{y-T_{\mathrm{S}}(G_{0}%
)}{S(G_{0})} \biggr)   \biggr)   =\delta.
\end{eqnarray*}
Then $\delta<1-C(G_{0})$ implies (\ref{M<1}) and from Theorem~\ref{Teore02}
follows that $T_{\mathrm{S}}$ and $S$ are weakly continuous. Since $S$ is
weekly continuous, Theorem~\ref{Teore01} implies that $T_{\mathrm{MM}}$ is
weakly continuous too, and so part (iii) follows.

Part (iv) follows from the formulas obtained in Section~\ref{secDerIF}.

(v) follows from  part (a)  of Theorem~\ref{Teorevarios} while part (vi)
follows from Lemma~\ref{lem1}. Part (vii) follows from (vi)  as was already
shown before stating the theorem.

To prove (viii) is enough to show that
\begin{equation}
n^{1/2}\mathrm{E}_{G_{n}}I_{T_{\mathrm{MM},\beta},G_{0}}(x,y)\rightarrow
_{d}N(0,V), \label{marmar}%
\end{equation}
where $V$ is given by (\ref{varerg}). From (\ref{INFMM1}), is immediate that
for all $\boldsymbol{\lambda}\in R^{q}, \lambda^{\prime}I_{T_{\mathrm{MM},\beta%
},G_{0}}(x_{i},y_{i})$ is a stationary ergodic martingale difference. Then
(\ref{marmar}) follows from the central limit theorem for martingale
differences (see, e.g., Theorem 23.1 of Billingsley~\cite{Billingsley1})
and the Cramer--Wald device.

Part (ix) will follow from
\begin{equation}
n^{1/2}\mathrm{E}_{G_{n}}I_{T_{\mathrm{MM}},G_{0}}(x,y)\rightarrow_{d}N(0,V),
\label{marmar1}%
\end{equation}
where $V$ is given by (\ref{varrmix}). According to (\ref{INFMM4}), we have
that
\[
I_{T_{\mathrm{MM}},G_{0}}(x_{i},y_{i} )= \frac{\sigma_{0}}%
{\mathrm{E}_{F_{0}}\psi_{1}^{\prime} (  (u-\alpha_{0})/\sigma_{0} )
}\psi_{1} \biggl(  \frac{u_{i}-\alpha_{0}}{\sigma_{0}} \biggr)  C_{0}%
^{-1}\underline{\dot{g}}(x,\beta_{0}),
\]
and therefore for all $\lambda\in R^{p+1},$ $\lambda^{\prime}%
I_{T_{\mathrm{MM}},G_{0}}(x_{i},y_{i} ) $ is a $\phi$-mixing process
 with mean 0  satisfying $%
{\sum_{i=1}^{\infty}}
\phi_{n}^{1/2}<\infty.$ Then by Theorem  20.1 of Billingsley
\cite{Billingsley1}, we have that $n^{1/2}\lambda^{\prime}\mathrm{E}_{G_{n}%
}I_{T_{\mathrm{MM}},G_{0}}(x, y)\rightarrow_{d}N(0,\lambda^{\prime}V\lambda),$
where%
\[
V=%
{\displaystyle\sum\limits_{i=-\infty}^{\infty}}
\mathrm{E} [  I_{T_{\mathrm{MM}},G_{0}}(x_{1},y_{1} ) %
I_{T_{\mathrm{MM}},G_{0}}^{\prime}(x_{1+i},y_{1+i} {)} ]  .
\]
Finally, the proof is completed noting that
\[
\mathrm{E} [  I_{T_{\mathrm{MM}},G_{0}}(x_{1},y_{1} {)}%
I_{T_{\mathrm{MM}},G_{0}}^{\prime}(x_{1+i},y_{1+i} {)} ]
=\frac{\sigma_{0}^{2} c_{i}}{\mathrm{E}_{F_{0}}^{2}\psi_{1}^{\prime
} (  (u-\alpha_{0})/\sigma_{0} )  }C_{0}^{-1}C_{i}C_{0}^{-1}%
\]
and using the Cramer--Wald device.

\subsubsection{\texorpdfstring{Proof of Theorem \protect\ref{TeoLtotal}}{Proof of Theorem 7}}

It is  completely similar to the proof of Theorem~\ref{teoremaso}. The only
differences  are that for part (iii) we use that in the case of a location
model we have $c (  G_{0} )  =0,$ and therefore condition (\ref{M<1})
reduces to $M_{G_{0}}(T_{\mathrm{M}}(G_{0}))<1.$ Note that this inequality is
implied by the condition  that $T_{\mathrm{M}}(G_{0})$ is well defined. So,
for this case, (\ref{M<1}) always holds, and that for part (iv) we use part
(b) of Theorem~\ref{Teorevarios} instead of part (a).

\section*{Acknowledgements}
This research was partially supported by Grants X--018  from
University of Buenos Aires,  PIP 112-200801-00216 and 0592  from CONICET and PICT  00899 and 0083
from ANPCyT.  We thank the Referees and Associate Editor for their comments
and suggestions which led to a much improved version of the paper.


\printhistory

\end{document}